\documentclass[11 pt]{amsart}

\usepackage{amssymb,color}
\usepackage{latexsym}
\usepackage{graphicx}
\usepackage{amsmath}
\usepackage{amsfonts}
\usepackage{natbib,pdfsync}

\allowdisplaybreaks
\iftrue 
\usepackage{amsmath}
\usepackage{amsfonts}
\newcommand{\field}[1]{\mathbb{#1}}
\DeclareMathOperator{\PR}{\field{P}}             
\DeclareMathOperator{\E}{\field{E}}              
\def\N{\field{N}}                                
\def\R{\field{R}}                                
\def\F{\field{F}}                                

\else
\def\PR{\mathop{\rm I\kern -0.20em P}\nolimits}  
\def\E{\mathop{\rm I\kern -0.20em E}\nolimits}   
\def\N{\mathop{\rm I\kern -0.20em N}\nolimits}   
\def\R{\mathop{\rm I\kern -0.20em R}\nolimits}   
\def\F{\mathop{\rm I\kern -0.20em F}\nolimits}   
\fi

\newtheorem{thm}{Theorem}[section]
\newtheorem{cor}[thm]{Corollary}
\newtheorem{lem}[thm]{Lemma}
\newtheorem{prop}[thm]{Proposition}
\newtheorem{ex}[thm]{Example}
\newtheorem{rem}[thm]{Remark}

\numberwithin{equation}{section}

\title[ Aggregation of Risks ]{Aggregation of rapidly varying risks 
 and asymptotic independence } 

\author[ A.\ Mitra ]{Abhimanyu Mitra}
\address{Abhimanyu Mitra\\School of OR\&IE, Cornell University,
Ithaca, NY-14853} \email{am492@cornell.edu}

\author[ S. I.\ Resnick ]{Sidney I. Resnick}
\address{Sidney I. Resnick\\School of OR\&IE, Cornell University,
Ithaca, NY-14853} \email{sir1@cornell.edu}

\thanks{S. I. Resnick and A. Mitra were partially supported by ARO
  Contract W911NF-07-1-0078  at Cornell University.}


\normalsize

\begin{document}

\begin{center}
\maketitle

\end{center}
\normalsize

\begin{abstract}

  We study the tail behavior of the distribution of the 
sum of asymptotically
  independent risks whose marginal distributions belong to the maximal
  domain of attraction of the Gumbel distribution. We impose
  conditions on the distribution of the risks $(X,Y)$ such that $P(X +
  Y > x) \sim (const)P (X > x)$. With the further assumption of
  non-negativity of the risks, the result is extended to more than two
  risks. We note a sufficient condition for a distribution 
to belong to both the maximal domain of attraction of
  the Gumbel distribution and the subexponential class. We provide
  examples of distributions which satisfy our assumptions.
The examples include cases where the marginal distributions
  of $X$ and $Y$ are subexponential and also cases where they are
  not. In addition, the asymptotic behavior of linear combinations of
  such risks with positive coefficients is explored leading to an
 approximate solution of an optimization problem which is applied
 to portfolio design.

\end{abstract}

\section{Introduction}

Estimating the probability that a sum of risks $X + Y$ exceeds a large
threshold is important in finance and insurance, and hence much
applied probability research has been dedicated to this goal. Recent
results are found in \citet{albrecher:asmussen:kortschak,
  kluppelberg:resnick, wang:tang:2006, asmussen:nandayapa,
  alink:lowe:wuthrich:2004, embrechts:puccetti:2006,
  ko:tang:2008}. Approximating this probability helps us evaluate
risk measures for investment portfolios as well as estimating
credit risk.

The problem is reasonably well understood when risks have regularly
varying marginal distributions but another important
large class of risk distributions
is the maximal domain of attraction of the
Gumbel distribution, denoted  $MDA(\Lambda)$, where
$$
\Lambda ( x) = \exp \{ - e^{-x} \}, \quad x \in \mathbb{R},
$$
 and $MDA(\Lambda)$ is the class of
distributions $F$ for which there exist $a_n>0, b_n \in \R$ such that
\begin{equation}\label{eqn:doaGumbel}
\lim_{n\to\infty}
n(1-F(a_nx +b_n)) =\lim_{n\to\infty}n\bar F(a_nx+b_n) =e^{-x},\quad
x\in\R
\end{equation}
\cite[page 38]{resnick:1987}. It is also well known that the risks having distribution in $MDA(\Lambda)$ are rapidly varying, i.e. $- \infty$-varying \cite[page 53]{resnick:1987}. Within the class of risks $(X, Y)$ with
marginal distributions $F, G \in MDA(\Lambda)$,  results on
aggregation of risks are known when $X$ and $Y$ are independent. 
However, actual risks are often
not independent and a somewhat weaker concept called
{\it asymptotic independence\/}, allows risks to be modeled as 
 dependent 
and  is more practical in many  modeling situations. Risks 
$X$ and $Y$ in a maximal domain of attraction are {\it asymptotically independent\/} if for all
${\bf{x}} = (x_1, x_2)$,
\begin{equation}\label{eqn:asyindep}
 \lim_{n \to \infty}H^n (a^{(1)}_n x_1 + b^{(1)}_n , a^{(2)}_n x_2 +
 b^{(2)}_n ) = G_1(x_1)G_2(x_2)
\end{equation}
  where  $H$ is the joint distribution of $X$ and $Y$ and both $G_1$
  and $G_2$ are non-degenerate extreme value distributions \cite[page
  229]{dehaan:ferreira:2006}. 
There are also results on aggregation of risks in the
 absence of  asymptotic
  independence where the analogue of \eqref{eqn:asyindep} holds but
  with a limit distribution which is not a product; see
 \cite{kluppelberg:resnick}. 
  
 This paper considers the case where the risks $X,Y$ are
 asymptotically independent with marginal distributions $F, G \in
 MDA(\Lambda)$.  We also allow one marginal tail to be lighter and the
 distribution with lighter tail does not necessarily belong to the
 maximal domain of attraction of the Gumbel distribution.

 Within the class of vectors $(X, Y)$ satisfying asymptotic
 independence and marginal distributions $F, G \in MDA(\Lambda)$, two
 prominent but very distinct behaviors have been observed.
\begin{enumerate}
\item First, suppose $(X, Y)$ are two iid risks with common
  distribution $F$ which is subexponential and $F \in
  MDA(\Lambda)$. Then $X$ and $Y$ are certainly asymptotically
  independent and
\begin{equation}\label{eqn:sharp2}
 \lim_{x \to \infty} \frac{P(X + Y > x)}{ P( X > x)} = 2. 
\end{equation}
So one possible behavior is that the sum has a distribution which is 
tail
equivalent to the distribution of a summand.

\item Very different tail behavior is exhibited 
in  Theorem 2.10 of \cite{albrecher:asmussen:kortschak}, who
exhibit a distribution of $(X, Y)$, with $X$ and $Y$ being
asymptotically independent and identically distributed with common
distribution $F \in MDA(\Lambda)$, but
\[ \lim_{ x \to \infty} \frac{ P( X + Y > x)}{P( X > x )} = \infty. \]
\end{enumerate}

In Section 2, we give a set of conditions on the joint distribution of
$(X, Y)$, guaranteeing behavior of the first sort, namely,
\begin{equation}\label{eqn:sharp}
 \lim_{x \to \infty} \frac{P(X + Y > x)}{ P( X > x)} =  1 + c  ,
\end{equation}
where  $ c = \lim_{ x \to \infty} {P( Y > x)}/{P( X > x)}$, the
limit being assumed to exist. If $c \in (0, \infty)$,
our conditions imply that $X,Y$ are
asymptotically independent and each belongs to the maximal domain of
attraction of the Gumbel. When $X,Y$ are identically distributed,
\eqref{eqn:sharp2} holds.
Under the further assumption of non-negativity of risks, the
result is extended for the case of more than two risks. In Section 3,
we provide examples of distributions which satisfy our conditions. The
examples include cases where the marginal distributions of 
$X$ and $Y$ are subexponential and also cases where they are not. We
also show one example which does not satisfy our  conditions 
but yet exhibits the tail equivalence between the distribution of the
sum and that of the summand.
Thus, our conditions are only sufficient. In Section 4, we
summarize asymptotic behavior of finite linear combinations of 
risks with non-negative coefficients. In Section 5, we suggest
approximate solutions for an optimization problem which is related to 
portfolio  design. The paper closes with concluding remarks
and a brief summary of numerical
experiments which give a feel for whether asymptotic equivalence is a
suitable numerical approximation for exceedance probabilities of
aggregated risks.

\section{Asymptotic tail probability for aggregated risk}

\subsection{Asymptotic tail probability for the sum of two random variables}
We give conditions guaranteeing \eqref{eqn:sharp}. The constant $c$
satisfies $c = \lim_{x\to\infty}{P( Y > x)}/{P( X > x)} \in [0,
\infty)$. When $c \in (0, \infty)$, $X$ and $Y$ are called
tail-equivalent \citep{resnick:1971} and then
 our conditions
guarantee that both the marginal distributions $F,G \in MDA(\Lambda)$
and $X$ and $Y$ are asymptotically independent. When $c = 0$, our
result extends to the case where
 $G$, the marginal distribution
of $Y$, does not belong to the maximal domain of attraction of the
Gumbel distribution and where $X$ and $Y$ need not be
asymptotically independent.

\subsubsection{Assumptions.} \label{ass_sec}
Suppose $(X,Y)$ is a pair of random variables satisfying the
following set of assumptions.
\begin{enumerate}
\item \label{Gumbel} The random variable $X$  has a distribution $F$
  whose right endpoint $x_0$ is infinite; that is,
\begin{equation}\label{rightend}
x_0 = \sup \{ x : F(x) < 1 \} = \infty. 
\end{equation}
Further $F
  \in MDA(\Lambda)$ so that \eqref{eqn:doaGumbel} is satisfied 
with centering constants $b_n\in \mathbb{R}$ and scaling
  constants $a_n>0$.
Equivalently (\citet{dehaan:1970}, \citet[page 28, 40-43]{resnick:1987} )
there exists a self-neglecting auxiliary function
$f(\cdot)$ with its derivative converging to 0, such that
\begin{equation}\label{eqn:selfneg}
  \lim_{ t \to \infty} \frac{\bar F (t + x f(t))}{\bar F (t)} =
  e^{-x}. 
\end{equation}
\item \label{ratioxy} The random variables $X$ and $Y$ have
distribution functions $F$ and $G$ such that
\[ \lim_{ x \rightarrow \infty}\frac{\bar G (x) } {\bar F (x) } = c  \in [0, \infty) .\]
\item \label{firstconx} The conditional distribution of $Y$ given $X >
  x$, satisfies for all $t>0$,
\[ \lim_{x \rightarrow \infty }P( |Y| > tf(x) | X > x ) = 0,
 \] 
where $f(x)$ is the auxiliary function corresponding to the
distribution of $X$ given in \eqref{eqn:selfneg},
\item \label{firstcony} and symmetrically assume for all $t>0$,
\[ \lim_{x \rightarrow \infty}P( |X| > tf(x) | Y > x ) = 0. \]
\item \label{secondcon} For some $L > 0$, suppose
\[ \lim_{x \rightarrow \infty}\frac{P( Y > Lf(x), X > Lf(x) )}{P(X >x)} = 0. \]
\end{enumerate}

\subsubsection{The main result.} The assumptions allow us to conclude
aggregated risks are essentially tail equivalent to individual risks.

\begin{thm}\label{thm:te}
 Under  Assumptions \ref{Gumbel}--\ref{secondcon} in
  Section \ref{ass_sec}, we have
\[P( X + Y > x) \sim (1 + c)P( X > x),\qquad x\to\infty.\]
\end{thm}

\subsubsection{Comments on the assumptions.}\label{subsubsec:comments}
Before giving a proof of Theorem \ref{thm:te},
we discuss implications of the assumptions.

\begin{rem} \label{rem:asy_ind}
\begin{enumerate}
\item \label{ab}
{\rm 
When $F \in MDA(\Lambda),$ we may choose $a_n, b_n$ appearing in
\eqref{eqn:doaGumbel}
as 
$b_n=b_F(n),$ $a_n=f(b_n)$. See \cite[page 40]{resnick:1987} or
\cite{dehaan:ferreira:2006}.

\item \label{asy_ind} 
If $c \in (0, \infty)$, then our assumptions guarantee both 
marginal distributions $F, G \in MDA(\Lambda)$ and also that $(X,Y)$ are
asymptotically independent. From Assumption \ref{Gumbel},  $F
\in MDA(\Lambda)$ and since $F$ and $G$ are tail-equivalent,
 from \cite{resnick:1971} we get that $G \in MDA(\Lambda)$. 
For asymptotic independence,
define, 
\begin{equation}\label{eqn:quantileFct}
b_F(t) = \inf\{ s : \frac{1}{1 - F} (s) \ge t
\}=\Bigl(\frac{1}{1-F}\Bigr)^\leftarrow (t),
\end{equation}
 and
similarly $b_G(t)$. From \cite[page 229]{dehaan:ferreira:2006},
if $F,G \in MDA(\Lambda)$ and 
\begin{align}\label{dehaan}
 \lim_{t \rightarrow \infty } \frac{P( X > b_F(t), Y > b_G(t) )}{P( X
   > b_F(t) )} &= 0, 
 \end{align}
then $(X,Y)$ are asymptotically independent according to \eqref{eqn:asyindep}.
When $c \in (0, \infty)$,  Assumption \ref{firstconx} implies
\eqref{dehaan}. To verify this, note first that 
Assumption \ref{firstconx} implies 
\begin{equation}\label{asymp_ind}
\lim_{x \rightarrow \infty } \frac{P( X > x, Y > x )}{P( X > x)} \le
\lim_{x \rightarrow \infty } \frac{P( X > f(x), Y > x )}{P( X > x)}=0,
\end{equation} since
$f(x)/x \to 0$ as $x\to\infty$ \citep[page 40]{resnick:1987}.
If  
 $c > 1$, then for sufficiently large $t$, $b_F(t) \le b_G(t)$ and
 therefore, using \eqref{asymp_ind},
\begin{align*}
\lim_{t \rightarrow \infty } \frac{P( X > b_F(t), Y > b_G(t) )}{P( X >
  b_F(t) )} &\le  \lim_{t \rightarrow \infty } \frac{P( X > b_F(t), Y
  > b_F(t) )}{P( X > b_F(t) )} \\ 
&=  \lim_{t \rightarrow \infty } \frac{P( X > t, Y > t )}{P( X > t)} =0,
\end{align*}
as required. A similar verification can be constructed for the case $0 < c < 1$.
For $c = 1, b_F(t) \sim b_G(t)$. Hence, 
\[ \frac{f(b_F(t))}{b_G(t)} \sim  \frac{f(b_F(t))}{b_F(t)}  \rightarrow 0.\]
So,
\begin{align*}
\lim_{t \rightarrow \infty } \frac{P( X > b_F(t), Y > b_G(t) )}{P( X > b_F(t) )} &\le  \lim_{t \rightarrow \infty } \frac{P( X > b_F(t), Y > f(b_F(t)) )}{P( X > b_F(t) )} \\
&= 0 \hskip 0.5 cm \hbox{ (by Assumption \ref{firstconx} and \eqref{rightend})}.
\end{align*}

\item \label{mean_excess}
 The auxiliary function $f(\cdot)$ can be replaced by any
 asymptotically equivalent function $\tilde f(\cdot)$; that is, if
 $
\lim_{x\to\infty}\tilde f(x) / f(x)=1,$ 
and if Assumptions \ref{firstconx}, \ref{firstcony}, \ref{secondcon} hold with
$f(\cdot)$, they also hold with $\tilde f(\cdot)$ replacing  $ f(\cdot)$ and
vice versa. Since the mean excess function
$$e (x) = E ( X - x | X > x)$$
 is asymptotically equivalent to any auxiliary function $f(x)$
(\cite[page 143]{embrechts:kluppelberg:mikosch:1997}, \cite[page
48]{resnick:1987}), $e(x)$ can also be taken as an auxiliary function.

\item \label{c=0}
 If $c = \lim_{ x
\rightarrow \infty}
{\bar G (x)}/{\bar F (x)} = 0,$
 we do not need Assumption \ref{firstcony} to conclude our result.

\item \label{notall_L}
An easier proof of the result can be given if Assumption \ref{secondcon} holds for all $L > 0$. But here we provide an example to show the importance of the weak version of Assumption \ref{secondcon}.
\begin{ex}
\[ X = -\log (U), \hskip 1 cm Y = -\log (1-U), \hskip 1 cm U \sim
 {\rm Uniform} \,(0,1) \]
{ \rm It is obvious that in this case both $X$ and $Y$ have
  distribution Exponential$(1)$. So, in this case, the auxiliary
  function is $f(x) =1$. Choose $L$ such that $\exp(-L) =
  \frac{3}{4}$, and
\begin{align*}
\frac{P( X > Lf(x), Y > L f(x) )}{P( X > x)} &= \frac{P( U < \exp(-L), 1-U < \exp(-L) )}{P( X > x)}\\
&= \frac{P( \frac{1}{4} < U <  \frac{3}{4} )}{P( X > x)}
= \frac{1}{2P( X > x)} \rightarrow \infty.
\end{align*}
Therefore, this particular choice of $L$ does not satisfy Assumption
\ref{secondcon}. 
The distribution of $(X, Y)$ is a special case of Example
\ref{heavylight} which discusses certain $L$ which do
 satisfy assumption \ref{secondcon}.
}
\end{ex}

\item \label{all_L}
 If, however, both $X$ and $Y$ are non-negative risks, and Assumption
 \ref{secondcon} is strengthened to hold for all $L > 0$,
then Assumptions
 \ref{firstconx} and \ref{firstcony} will be automatically
 satisfied. The proof of
 this  follows from  $\lim_{x \to \infty} f(x)/x
 =0.$ 

\item \label{comparison}
Similar limit results are found in Lemma 2.7 of
\cite{albrecher:asmussen:kortschak}  and Theorem 2.1 of
\cite{ko:tang:2008}. They have assumed that one of the marginal
distributions of the two asymptotically independent variables $X$ and
$Y$, say the distribution of $X$, is subexponential, (i.e. $X \in
\mathcal{S}$, where $\mathcal{S}$ is the set of all subexponential
distributions) and worked on finding conditions for the
tail-equivalence of the marginal distribution of $X$ and the sum $X +
Y$. 
Our assumptions are different:  We assume that one of the marginal
distributions of the two asymptotically independent variables $X$ and
$Y$, say the distribution of $X$, belongs to the domain of attraction
of Gumbel, i.e. $X \in MDA( \Lambda)$. We do not assume the marginal
distribution of $X$ is subexponential.  

In examples where the marginal distributions of the two asymptotically
independent and identically distributed random variables $X$ and $Y$
belong to the class $MDA(\Lambda) \cap \mathcal{S}$,
an issue is the relative strength of our conditions versus those of
Theorem 2.1 of \cite{ko:tang:2008} .  We can not show either
set of conditions implies the other. Below we present an example which
satisfies
 our set of conditions, but does not satisfy the set of conditions
 given in Theorem 2.1 of \cite{ko:tang:2008}.  Thus our set of
 conditions is not stronger. 

\begin{ex}
{\rm Suppose,
$ X = \exp(X_1) ,$ $ Y = \exp(X_2), $
where $(X_1, X_2)$ is bivariate normal with correlation $\rho \in (0,1)$.
For simplicity, assume  each $X_i$ has mean 0 and variance $1$. It is well known that lognormal distribution belongs to the class $MDA(\Lambda) \cap \mathcal{S}$. In Example \ref{lognormal}, we show $(X, Y)$ satisfy our set of conditions. Here we show that this example does not satisfy Assumption 2.1 of \cite{ko:tang:2008}, i.e. for all $x^* > 0$,
\begin{equation}\label{kotangcondn}
\limsup_{x \to \infty} \sup_{ x^* \le t \le x} \frac{P( Y > x - t | X = t )}{ P( Y > x - t )} = \infty.
\end{equation}
From the exchangeability of $X$ and $Y$, it is obvious that \eqref{kotangcondn} holds even if the role of $X$ and $Y$ is interchanged. 

\begin{align}\label{kotangcondn2}
\sup_{ x^* \le t \le x}  \frac{P( Y > x - t | X = t )}{ P( Y > x - t )} = \sup_{ x^* \le t \le x} \frac{ \bar \Phi \left( \frac{ \log(x - t) - \rho \log t }{ \sqrt{ 1- \rho^2}} \right)}{ \bar \Phi \left( \log (x-t) \right)} 
& \ge \frac{ \bar \Phi \left( \frac{ \log(x/2) - \rho \log (x/2) }{ \sqrt{ 1- \rho^2}} \right)}{ \bar \Phi \left( \log (x/2) \right)} \nonumber \\
&= \frac{ \bar \Phi \left( \frac{ 1 - \rho }{ \sqrt{ 1- \rho^2}} \log(x/2) \right)}{ \bar \Phi \left( \log (x/2) \right)}\rightarrow \infty.
\end{align}
The inequality above follows from choosing $x$ enough large so that $x/2 > x^*$ and putting $t = x/2$. The last convergence follows from the fact that normal distribution belongs to the class $MDA(\Lambda)$ and hence $\bar \Phi$ is $-\infty$-varying \cite[page 53]{resnick:1987}. Note, $0 < \rho <1$ entails 
$ \frac{ 1 - \rho }{ \sqrt{ 1- \rho^2}} < 1$. Hence, from \eqref{kotangcondn2} it is obvious that \eqref{kotangcondn} holds.
}
\end{ex}
}
\end{enumerate}
\end{rem}

\subsubsection{Proof of Theorem \ref{thm:te}.}
We prove Theorem \ref{thm:te} using a Proposition and a
Lemma, which we prove first. Note, we do not need the assumption that the marginal distributions are sub-exponential, which is a necessary condition in the case where $X$ and $Y$ are independent.

\begin{prop} \label{prop1}
Under Assumptions 
\ref{Gumbel} and \ref{firstconx} of Section \ref{ass_sec}, 
we have
\begin{equation} \label{prop1eq1}
\lim_{ n \rightarrow \infty} P ( Y \le a_n z | X > a_n x + b_n ) =
    1_{\{z > 0\}},\qquad  z \ne 0 ,  \; x \in \R.
    \end{equation}
    and from Assumptions \ref{Gumbel} and \ref{firstcony} of Section
  \ref{ass_sec}, we have
\begin{equation}\label{eqn:prop2}
\lim_{ n \rightarrow \infty} P ( X \le a_n z | Y > a_n x + b_n ) =
    1_{\{z > 0\}},\qquad z \ne 0 ,  \, x \in \R.
\end{equation}

\end{prop}

\begin{proof}
  The self-neglecting property of the auxiliary function $f$, i.e. 
  \begin{equation*}
  \lim_{ t \to \infty} \frac{ f( t + x f(t))}{f(t)} = 1, \hskip 1 cm x \in \R,
  \end{equation*}
   implies that
  \begin{equation}\label{prop1:eq1}
  \lim_{ t \to \infty} P( |Y| > z f(t) | X > t) = \lim_{ t \to \infty} P( |Y| > z f(t) | X > t + x f(t) ). 
  \end{equation}
  Hence, by noting that $a_n = f(b_n)$ and $\lim_{n \to \infty} b_n = \infty$, the result follows from \eqref{prop1:eq1}. The second part is proved similarly.
\end{proof}

 \begin{lem}\label{lem:vague} 
 (i) Assumptions \ref{Gumbel}, \ref{ratioxy}, and \ref{firstconx} of Section \ref{ass_sec} imply that the sequence of measures 
\[ n P[a_n^{-1}( X - b_n, Y) \in (dx, dy) ] 
\] 
converges vaguely on $([ -M, \infty] \times [-\infty, \infty])$ as $n \to
\infty$, 
to the limit
measure 
$m_{1, \infty}( dx, dy) =  e^{-x}dx\epsilon_0(dy),$
for some $M >  L$(from Assumption \ref{secondcon} of Section \ref{ass_sec}) such that $-M$ is a continuity point of $\frac{X-b_n}{a_n}$ for all $n$.\\
(ii) Assumptions \ref{Gumbel}, \ref{ratioxy}, and \ref{firstcony} of Section \ref{ass_sec} imply that the sequence of measures 
\[ n P[a_n^{-1}( Y - b_n, X) \in (dx, dy) ] 
\] 
converges vaguely on $([ -M, \infty] \times [-\infty, \infty])$ as $n \to
\infty$, 
to the limit
measure 
$m_{2, \infty}( dx, dy) =  ce^{-x}dx\epsilon_0(dy),$
for some $M >  L$(from Assumption \ref{secondcon} of Section \ref{ass_sec}) such that $-M$ is a continuity point of $\frac{Y-b_n}{a_n}$ for all $n$.

\end{lem}

\begin{rem}\label{sameM}
 \rm{Since all the discontinuity points of $\frac{X - b_n}{a_n}$ for all $n$ is countable, choice of such an $M > L$ is not a problem. Moreover, the $M$ in the two parts of the lemma (i) and (ii) may be chosen to be the same.}
\end{rem}
\begin{proof}  We consider convergence of the measures evaluated on
  certain relatively compact regions which
guarantee vague convergence. 

\smallskip
{\sc Region 1:} $(x, \infty] \times [-\infty, y]$,  $x > -M, y \ne 0$. 
As $n \to \infty$,
\begin{align*}
nP\Bigl[ \frac{X -b_n}{a_n} >x, \frac{Y}{a_n} \le y \Bigr]
&= nP\left[\frac{X -b_n}{a_n}> x \right] P\left[ \frac{Y}{a_n} \le y
  \big | \frac{X -b_n}{a_n} > x \right] \\ 
&\rightarrow e^{-x}1_{\{y>0 \}} 
= m_{1, \infty}( (x, \infty] \times [-\infty, y] )
\end{align*}
by Proposition \ref{prop1}.

{\sc Region 2:} $[-M, x] \times (y, \infty]$, $x > -M, y \ne 0$.
Since $-M$ is a continuity point of $\frac{X-b_n}{a_n}$ for all $n$, as $n \to \infty$, 
\begin{align*}
nP\Bigl[-M \le \frac{X -b_n}{a_n} \le x, \frac{Y}{a_n} > y \Bigr] &= nP\Bigl[-M < \frac{X -b_n}{a_n} \le x, \frac{Y}{a_n} > y \Bigr] \\
&= nP\Bigl[ \frac{X -b_n}{a_n} > -M, \frac{Y}{a_n} > y \Bigr] - nP\Bigl[ \frac{X -b_n}{a_n} > x, \frac{Y}{a_n} > y \Bigr]\\
&=  nP\left[\frac{X -b_n}{a_n}> -M \right] P\left[ \frac{Y}{a_n} > y
  \big | \frac{X -b_n}{a_n} > -M \right] \\
& \qquad -  nP\left[\frac{X -b_n}{a_n}> x \right] P\left[ \frac{Y}{a_n} > y
  \big | \frac{X -b_n}{a_n} > x \right] \\
& \qquad \rightarrow (e^M -e^{-x}) 1_{ \{y < 0\} } = m_{1, \infty}( [-M, x] \times (y, \infty] ),
\end{align*}
by Assumption \ref{Gumbel} and Proposition \ref{prop1}. 

Arguments for convergence on the following regions follow in a
similar fashion using Proposition \ref{prop1} :\\
\indent {\sc Region 3:} $(x, \infty] \times (y, \infty] $, \hskip 0.1 cm $x > -M, y \ne 0$,\\
\indent {\sc Region 4:} $[-M, x] \times [-\infty, y] $, \hskip 0.1 cm $x > -M, y \ne 0$.\\
This concludes the proof of vague convergence on part (i).

The proof of part (ii) is similar, only notice that if c = 0, we do not need the Assumption \ref{firstcony}. In this case, note that the limit measure $m_{2, \infty}( dx, dy)$ is a zero measure. 
Also note, using Assumptions \ref{Gumbel} and \ref{ratioxy}, we get
\begin{align*}
nP\left[\frac{Y -b_n}{a_n} \ge -M \right] \rightarrow ce^M = 0,
\end{align*}
which is enough to prove the convergence in this case.
\end{proof}

This leads to a formal statement of the main result.
\begin{thm} \label{mainthm}
Under the Assumptions in Section \ref{ass_sec},
 \begin{equation}\label{thmeqn3}
 \lim_{ x \rightarrow \infty} \frac{P(X +Y >x)}{P(X > x)} = (1+ c).
 \end{equation}
\end{thm}

\begin{proof} Choose $M$ as in Remark \ref{sameM}. We split $P( X+Y > b_n)$ as 
\begin{align}\label{mainthm1}
P(X + Y > b_n) &= P( X + Y > b_n, X > b_n - Ma_n ) + P( X + Y > b_n, Y > b_n - Ma_n)  \nonumber\\
& \quad - P( X + Y > b_n, X > b_n - Ma_n, Y > b_n - Ma_n) \nonumber \\
& \qquad + P( X + Y > b_n, X \le b_n - Ma_n, Y \le b_n - Ma_n).
\end{align}
Using Assumption \ref{Gumbel} and \eqref{asymp_ind}, we get
\begin{align}\label{mainthm2}
nP( X + Y > b_n, X > b_n - Ma_n, Y > b_n - Ma_n) \le nP( X > b_n - Ma_n, Y > b_n - Ma_n) \nonumber \\
= nP( X > b_n - Ma_n)\frac{P( X > b_n - Ma_n, Y > b_n - Ma_n)}{P( X > b_n - Ma_n)}
\rightarrow e^M . 0 = 0,
\end{align}
since $b_n - M a_n \to \infty$.
Now, consider the convergence of the last term of \eqref{mainthm1} mutiplied by n.
\begin{align}\label{mainthm3}
P( X + Y > b_n, X \le b_n - Ma_n, Y \le b_n - Ma_n) \le nP( X >  Ma_n, Y > Ma_n) \nonumber \\
\sim \frac{P( X > Mf(b_n), Y > Mf(b_n))}{P( X > b_n)} \le \frac{P( X > Lf(b_n), Y > Lf(b_n))}{P( X > b_n)} 
\rightarrow  0,
\end{align}

by \eqref{rightend} and Assumption \ref{secondcon}.

To deal with the first term of \eqref{mainthm1} mutiplied by n, we first define a function $T$ as T : $[ -M, \infty] \times [-\infty, \infty] \mapsto
 (\infty,\infty]$ by 
 $$T(x,y ) =
\begin{cases}
x+y,& \text{ if } y>-\infty,\\
0,& \text{ if } y=-\infty,
\end{cases}
$$
and hence
\begin{align}\label{mainthm4}
nP( X + Y > b_n, X > b_n - Ma_n) = nP( a_n^{-1}(X - b_n, Y) \in T^{\leftarrow}( (0, \infty]) \cap [(-M, \infty] \times \{0 \}] ).
\end{align}

Note, that every set in the space $[-M, \infty] \times [-\infty, \infty]$ is relatively compact, and hence so is $T^{\leftarrow}( (0, \infty]) \cap [(-M, \infty] \times \{0 \}] = S$ (say).
Also, since the limit measure $m_{1, \infty}$ is concentrated on $[-M, \infty] \times \{ 0\}$,
\begin{align}\label{mainthm5}
m_{1, \infty}( \delta S ) = m_{1, \infty}( \delta S \cap [ [-M, \infty] \times \{ 0\} ] ) = m_{1, \infty}( \{0\} \times \{ 0\}) = 0.
\end{align}
Hence, using Lemma \ref{lem:vague}, \eqref{mainthm4} and \eqref{mainthm5}, we get
\begin{align}\label{mainthm6}
nP( X + Y > b_n, X > b_n - Ma_n) \rightarrow m_{ 1, \infty}( S ) = 1.
\end{align}
Similarly,
\begin{align}\label{mainthm7}
nP( X + Y > b_n, Y > b_n - Ma_n) \rightarrow m_{ 2, \infty}( S ) = c.
\end{align}
Hence, using \eqref{mainthm1}, \eqref{mainthm2}, \eqref{mainthm3}, \eqref{mainthm6} and  \eqref{mainthm7}, we get, 
\begin{align*}
\lim_{ x \to \infty} \frac{P( X + Y > x)}{P( X > x)} = \lim_{n \to \infty} nP( X + Y > b_n) = 1 + c, 
\end{align*}
and we conclude our result. 
\end{proof}

One immediate application of Theorem \ref{mainthm}
 is to 
the subexponential family of distributions denoted $\mathcal{S}$. 
The class  $MDA(\Lambda) \cap \mathcal{S}$ has been studied
in \citep[page
149]{embrechts:kluppelberg:mikosch:1997} and several
sufficient conditions for belonging to this class are given 
in \cite{ goldie:resnick1988b}.  Corollary \ref{subexponential}
 gives an
additional sufficient condition and follows directly from Theorem \ref{mainthm}.
Example \ref{sub} exhibits a distribution which satisfies
the conditions of this Corollary. 
 
\begin{cor} \label{subexponential} Suppose, $F \in MDA(\Lambda)$ with
  auxiliary function $f(x)$ as described in Assumption \ref{Gumbel} of
Section  \ref{ass_sec}. Suppose, also, $\lim_{x \to \infty} f(x) =
\infty$, and for some $L > 0$,  
 \begin{equation} \label{subeqn1}
  \lim_{x \to \infty} \frac{ {\left[\bar F ( L f(x))\right]}^2} {\bar
    F (x)} = 0.  
  \end{equation}
 Then, for $X$ and $Y$ iid with common distribution $F$ we have, as $x \to \infty$,
 \begin{equation*}
 P[ X + Y > x ] \sim  2P[ X > x ], 
 \end{equation*}
 and therefore, if $F$ concentrates on $[0, \infty)$, $F \in \mathcal{S}$. \end{cor}

 Following Remark \ref{rem:asy_ind}(\ref{mean_excess}), it is enough
 to check \eqref{subeqn1}  with any  $\tilde 
 f(x)$  satisfying $ \tilde f(x) \sim f(x)$.  Note also it is natural
 to add the assumption $f(x)\to\infty$, since if $F\in
 MDA(\Lambda)\cap \mathcal{S}$, then necessarily $f(x)\to\infty$
 \citep{goldie:resnick1988b}. 
  
\subsection{Asymptotic tail probability for the sum of more than
two non-negative random variables}
Suppose, among the risks
$X_1, X_2, \dots X_d$, there is no heavier tail than $X_1$ in the
sense that it is not true that 
$$\lim_{x\to\infty} \frac{\bar F_i(x)}{\bar F_1(x)}=\infty,\quad
i=2,\dots,d.$$
Assume
  $X_1$ satisfies Assumption \ref{Gumbel}
 of Section \ref{ass_sec} and that $X_1, X_2, \dots
X_d$ pairwise satisfy the  Assumptions 
\ref{firstconx} and \ref{firstcony} of Section \ref{ass_sec} with the
auxiliary 
function $f(\cdot)$ of $X_1$. By this, we mean for all pairs $1 \le i
\ne j \le d$, and for $t > 0$, 
$$
\lim_{ x \to \infty} \frac{P( X_j > t f(x) , X_i > x )}{P(X_i
> x)} = 0,
$$
which implies
\begin{equation}\label{pairwise1}
\lim_{ x \to \infty} \frac{P( X_j > t f(x) , X_i > x)}{P(X_1
> x)} = 0.
\end{equation}
 Also, suppose, the risks $X_1, X_2, \ldots X_d$ pairwise satisfy 
Assumption \ref{secondcon} of Section \ref{ass_sec} with auxiliary function
 $f(\cdot)$ of $X_1$ so that 
 for $ 1 \le i < j \le d$, there exists some $L_{ij} > 0$, such that either
 \[ \lim_{x \to \infty} \frac{P( X_i > L_{ij} f(x), X_j > L_{ij} f(x) )}{P(X_i > x)} = 0,
 \]
 or, 
 \[ \lim_{x \to \infty} \frac{P( X_i > L_{ij} f(x), X_j > L_{ij} f(x) )}{P(X_j > x)} = 0.  \]
In either case, we have, for $1 \le i < j \le d$, for some $L_{ij} > 0$,
\begin{equation}\label{pairwise2}
\lim_{x \to \infty} \frac{P( X_i > L_{ij} f(x), X_j > L_{ij} f(x)
)}{P(X_1 > x)} = 0.  
\end{equation}
Under the additional assumption of non-negativity,
Theorem \ref{mainthm} can be extended to more than two risks.

\begin{cor} \label{sumcor}
Assume, $X_1, X_2, \ldots X_d$ are
non-negative random variables which pairwise satisfy Assumptions
\ref{firstconx}, \ref{firstcony}, \ref{secondcon} of Section
 \ref{ass_sec} with the auxiliary function $f(\cdot)$ of $X_1$. Moreover,
the distribution of $X_1$ satisfies Assumption \ref{Gumbel}
of Section \ref{ass_sec} and suppose
\begin{equation} \label{cor1eqn1}
\lim_ {x \rightarrow \infty} \frac{ P( X_i > x )}{ P( X_1 > x)} =
c_i \in [0, \infty), \qquad i = 2, 3, \ldots,d.
\end{equation}
Define, $S_j = X_1 + X_2 + \ldots X_j,  1\le j \le d$ and we have, for $x \in \R$,
\begin{equation}\label{eqn:d}
 \lim_{n \rightarrow \infty} n P (S_d > a_nx + b_n) = ( 1 + \sum_{i=2}^d c_i )e^{-x}
 \end{equation} and hence
 \begin{equation}\label{eqn:mored}
 \lim_{ x \rightarrow \infty} \frac{P(S_d >x)}{P(X_1 > x)} = ( 1 + \sum_{i=2}^d c_i )
 \end{equation}
\end{cor}

\begin{rem}\label{asyindpair} 
\begin{enumerate}
\item \label{item:asyindpair} {\rm
{\it Asymptotic independence of the random variables\/}:
 Suppose, for all $i$, $c_i \in (0, \infty)$. Then for any $1 \le i \ne j \le d$, the
 pair $(X_i,X_j) $  is asmptotically independent by 
Remark \ref{rem:asy_ind}(\ref{asy_ind}). Since the random variables are pairwise
asymptotically independent, they are also
 asymptotically independent \citep[page 291]{resnick:1987}.

\item \label{item:nonneg}
{\it Non-negativity of random variables\/}: The only additional assumption
added to the list in Section \ref{ass_sec} is that the random variables are non-negative.

\item \label{relax}
{\it Relaxation\/}: We have shown in \eqref{pairwise1} and \eqref{pairwise2} that pairwise satisfaction of Assumptions \ref{firstconx},
\ref{firstcony}, \ref{secondcon} of Section \ref{ass_sec}  implies
that for $ 1
\le i \ne j \le d$, for $t >  0$, 
\begin{eqnarray*}
\lim_{ x \to \infty} \frac{P( X_j > t f(x) , X_i > x)}{P(X_1
> x)} &=& 0,
\end{eqnarray*}
and for $ 1 \le i < j \le d$, there exists $L_{ij} > 0$,
\begin{eqnarray*}
\lim_{ x \to \infty} \frac{P( X_j > L_{ij} f(x) , X_i > L_{ij}
f(x))}{P(X_1
> x)} &=& 0.
\end{eqnarray*}
We will show that actually these conditions are enough to
get the desired conclusion. } 
\end{enumerate}\end{rem}

 \begin{proof}
 We prove the result by induction under the relaxation
Remark \ref{asyindpair}(\ref{relax}). The base case
of the induction for $d = 2$ is already proved in Theorem
\ref{mainthm}, so 
suppose, the result is true for $d = k \ge 2$ and we have
\begin{eqnarray}
 \lim_{n \rightarrow \infty} n P (S_k > a_nx + b_n) &=& ( 1 + \sum_{i=2}^k c_i )e^{-x}
 \end{eqnarray}
and
 \begin{eqnarray} \label{cor1eqn2}
 \lim_{x \rightarrow \infty} \frac{P( S_k > x)}{P( X_1 > x )} &=& 1
 + \sum_{i=2}^k c_i
 \end{eqnarray}
 Therefore, we have
  \begin{eqnarray} \label{sumcorratio}
 \lim_{x \rightarrow \infty} \frac{P( X_{k +1} > x)}{P( S_k > x )} &=& \frac{c_{k+1}}{1
 + \sum_{i=2}^k c_i} \in [0, \infty).
 \end{eqnarray}
We will use  Theorem \ref{mainthm} with $X =
S_k$ and $Y = X_{k +1}$. It remains to check the Assumptions in  Theorem \ref{mainthm}. For Assumption \ref{Gumbel}, note that $S_k$ is tail-equivalent to $X_1$ and use the fact that $F \in MDA(\Lambda)$ is closed under tail-equivalence. Assumption \ref{ratioxy} is already checked in \eqref{sumcorratio}.

Note that from the induction hypothesis $P( S_k > x) \sim P( \cup_{i=1}^k ( X_i > x))$ and from the positivity of the risks $[ S_k > x ] \supseteq  \cup_{i=1}^k [X_i > x] $. From these two facts, it easily follows that
\begin{equation}\label{eqn:2.27.5}
 \lim_{x \rightarrow \infty} \frac{P( (S_k > x) \cap
   {(\cup_{i=1}^k(X_i > x))}^c)}{P(S_k > x)} = 0.
\end{equation}

Since $S_k$ and $X_1$ are tail equivalent, by \cite{resnick:1971}, the auxiliary function $\tilde f(\cdot)$
of $S_k$ is asymptotically equal to  the
auxiliary function $f(\cdot)$ of $X_1$.
Therefore, given $\epsilon \in (0,1)$, there exists $T$ 
such that for all $x > T , \tilde f(x)
> \epsilon f(x).$ 
We now check Assumption \ref{firstconx}.
For $t  >  0, x > T$, using \eqref{eqn:2.27.5}, as $x\to\infty$,
 \begin{align*}
P( & |X_{k + 1}| > t \tilde f(x) | S_k > x ) \le P( X_{k + 1} > t
\epsilon f(x) | S_k > x ) \\ 
   &= \frac{P( X_{ k + 1 } > t \epsilon f(x), S_k > x)}{ P( S_k > x )} 
   \sim  \frac{P( X_{ k + 1 } > t \epsilon f(x), S_k > x, \cup_{i=1}^k
     \{ X_i > x\} )}{ P( S_k > x )}\\ 
   &\le \frac{P( X_{ k + 1 } > t \epsilon f(x), \cup_{i=1}^k \{ X_i >
     x\} )}{ P( S_k > x )}
\le \frac{\sum_{i = 1}^k P( X_{ k + 1 } > t \epsilon f(x), X_i > x )}{
  (1 + \sum_{i=2}^k c_i) P( X_1 > x )} 
    \rightarrow 0
   \end{align*}
   by \eqref{pairwise1}.

For Assumption \ref{firstcony}, if $c_{k+1} = 0,$ following Remark
\ref{rem:asy_ind}(\ref{c=0}), there is no need to check assumption
\ref{firstcony}. So, suppose, $c_{k+1} > 0.$ Then for any $t  >  0$,
 as $x\to\infty$,
\begin{align*} 
   P( |S_k| > t \tilde f(x) | &X_{ k + 1} > x ) \le P( S_k > t \epsilon f(x) | X_{ k + 1} > x )\\
   &\le \sum_{ i = 1}^{k}  P( X_i  > t \epsilon f(x)/k | X_{ k + 1} > x )\\
   &= \sum_{ i = 1}^{k}  \frac{P( X_i  > t \epsilon f(x)/k, X_{ k + 1}
     > x )}{P( X_1 > x )}\frac{P ( X_1 > x)}{P( X_{ k + 1} > x )} 
\rightarrow 0.
   \end{align*}
 For Assumption \ref{secondcon},
we know from the assumptions in the statement of Corollary
\ref{sumcor} that the random 
 variables satisfy Assumption \ref{secondcon} of Section
\ref{ass_sec} pairwise with auxiliary
 function $f(\cdot)$ of $X_1$. Thus, for $1 \le i < j \le d$,  \eqref{pairwise2} holds.
 We check Assumption \ref{secondcon} with  $L = kL_{max}/ \epsilon$, where $L_{max} = \max_{1 \le i \le k}L_{i, k+1}$ (recall, equation \eqref{pairwise2}).
 Then, for sufficiently large $x$, using $\tilde f(\cdot)$ as the auxiliary function of $S_k$,
  \begin{align*}
&   \frac{P( X_{ k +1} > L \tilde f(x) , S_k > L \tilde f(x) )}{ P(S_k > x )}
   \le \frac{P( X_{ k +1} > L \epsilon f(x) , S_k > L \epsilon f(x) )}{ P(S_k > x
   )}\\
   &\qquad \le \frac{P( X_{ k + 1 } > k L_{max} f(x), \cup_{i=1}^k \{ X_i > L_{max} f(x) \})}{ P( S_k > x )}\\
    &\qquad \le \frac{P( X_{ k + 1 } > L_{i, k + 1} f(x), \cup_{i=1}^k \{ X_i > L_{i, k + 1} f(x) \})}{ P( S_k > x )}\\
      &\qquad \le \frac{\sum_{i=1}^k P( X_{ k + 1 } > L_{i, k + 1}f(x), X_i > L_{i, k + 1}f(x) )}{ P( S_k > x )} \\
      &\qquad \sim \frac{\sum_{i=1}^k P( X_{ k + 1 } > L_{i, k + 1}f(x), X_i > L_{i, k + 1}f(x) )}{ (1 + \sum_{i=2}^k c_i )P( X_1 > x )} 
\rightarrow 0
   \end{align*}
 by \eqref{pairwise2}. This completes the induction proof.
\end{proof}

\section{Examples}
This section shows a few of the many models that satisfy the
Assumptions in Section \ref{ass_sec}.  In all examples, both $X$ and $Y$ are
non-negative random variables and it is straightforward to extend
these examples to the
$d$-dimensional case and show the assumptions of Corollary
\ref{sumcor} are satisfied. 

Our conditions are only sufficient and we exhibit one example where
our conditions do not hold, but tail equivalence in Theorem
\ref{mainthm} holds true. Finding a necessary and sufficient condition
for the conclusion of Theorem \ref{mainthm} is still an open but subtle and difficult issue.

\begin{ex}  \label{minxy} 
{\rm Suppose  $X_1, X_2, X_3 $ are iid  with common distribution $F$, where for $\alpha > 1$,
\[ \bar F(x) = \left\{ \begin{array}{cc}
\exp \{ - {(\log x)}^{\alpha} \}, & \hbox{if} \hskip 0.2 cm x > 1, \\
1, & \hbox{if} \hskip 0.2 cm x \le 1. \end{array} \right.\]
Define,
\[ X = X_1 \wedge X_2, \hskip 1 cm Y= X_2 \wedge X_3 \]
It is easy to check $X$ and $Y$ are identically distributed with the
common distribution  $F_1$, where
\[ \bar F_1 (x) = \exp ( - 2{(\log x)}^{\alpha}), \hskip 0.2 cm x > 1. \\
\]
It can be checked that $F_1$ is a Von-Mises function; that is, it satisfies,
$$\frac{\bar F_1 F_1^{''}}{(F_1^{'})^2}\to -1,$$
a sufficient condition for $F_1\in MDA(\Lambda)$,
and 
\[ f(x) = \frac{\bar F_1 (x)}{  F_1^{'} (x)}=
\frac{x}{2 \alpha {(\log x)}^{ \alpha -1} }, \hskip 1 cm x > 1\]
serves as an auxiliary function \citep[page 40]{resnick:1987}. Also, \eqref{rightend} is obvious and therefore, Assumption \ref{Gumbel} of Section 
\ref{ass_sec} is satisfied.
Checking Assumption \ref{ratioxy} is straightforward, so consider
Asumption \ref{firstconx}. Fix $t > 0$, recall $f(x)/x\to 0$, and note as $x \to \infty$,
\begin{align*}
&\frac{P( X > x, Y > t f(x))}{ P(X > x)} = \frac{P( X_1 > x, X_2
> x \vee t f(x), X_3 >  t f (x))}{ P(X_1 > x, X_2 > x)}\\
&\qquad \sim \frac{P( X_1 > x, X_2 > x, X_3> t f(x))}{ P(X_1 >  x, X_2 > x)} 
= P(X_3 > t f(x))\rightarrow 0,
\end{align*}
since $f(x) \to \infty$. Assumption \ref{firstcony} is verified the
same way. For Assumption \ref{secondcon}, we have with $L = 1$, 
\begin{align} \label{exeqn1}
&\frac{P( X > f(x), Y > f(x))}{ P(X > x)} = \frac{P( X_1 >
f(x), X_2 > f(x), X_3 > f(x))}{ P(X_1 > x, X_2 > x)} \nonumber\\
 &\qquad = \frac{ {\bar F (f(x))}^3}{{\bar F (x)}^2} 
= \exp\left\{ - \left[ 3{(\log f (x) )}^{\alpha} - 2{(\log x )}^{\alpha} \right]\right\} \nonumber \\
&\qquad = \exp\left\{ - 2{(\log x )}^{\alpha} \left[ \frac{3}{2}{\left(\frac{ \log f (x) }{\log x }\right) }^{\alpha}- 1\right]\right\} \nonumber \\
&\qquad = \exp\left\{ - 2{(\log x )}^{\alpha} \left[ \frac{3}{2}{\left( 1 - \frac{ \log( 2 \alpha {(\log x)}^{ \alpha -1} )}{\log x }\right) }^{\alpha}- 1\right]\right\} .
\end{align}
Since the exponent in \eqref{exeqn1} converges to $- \infty$ as $x \to
\infty$, Assumption \ref{secondcon} is satisfied
and this pair $(X, Y)$ satisfies the Assumptions in Section
\ref{ass_sec}. 
}
\end{ex} 

\begin{ex} \label{sub} 
{\rm Suppose, $X$ and $Y$ are independent and identically distributed  with common distribution $F$, where for $ \alpha > 1$,
\[ \bar F(x) = \left\{ \begin{array}{cc}
\exp ( - {(\log x)}^{\alpha}) & \hbox{if} \hskip 0.2 cm x > 1, \\
1 & \hbox{if} \hskip 0.2 cm x \le 1. \end{array} \right. \]
As in Example \ref{minxy}, one can check the subexponentiality condition
\eqref{subeqn1}
 with $L = 1$ and by Corollary \ref{subexponential},  $F$ is
 subexponential. Hence, 
 \[ P( X + Y > x ) \sim 2P( X > x ). \]  
}
\end{ex}

\begin{ex} \label{rhoextreme} 
{\rm Suppose, $X \sim \rm{Lognormal} ( \mu , \sigma^2 )$
and $Y = {e^{2\mu}}/{X}$ so that $X \stackrel{d}{=} Y$. 
We check  the Assumptions in Section
\ref{ass_sec} for the pair $(X, Y)$. 
The distribution $\rm{Lognormal} ( \mu , \sigma^2 )$ belongs to the maximal
domain of attraction of the Gumbel distribution and its mean excess
function $e(x)$ has 
the form \citep[page 147, 161]{embrechts:kluppelberg:mikosch:1997} 
\[ e(x) = \frac{\sigma^2 x}{ \log x - \mu} ( 1 + o(1)).\]
Also, \eqref{rightend} is obvious and so, Assumption \ref{Gumbel} of Section \ref{ass_sec} is true.
Following Remark \ref{rem:asy_ind}(\ref{mean_excess}) and the
 form of $e(x)$, we may assume the auxiliary function
\[ f(x) = \frac{\sigma^2 x}{ \log x - \mu}. \] 
To verify Assumption \ref{firstconx}, fix  $t > 0$, and note as $x \to \infty$, 
$$
\frac{P( X > x, Y > t f(x))}{ P(X > x)} = \frac{P( X  > x, e^{2\mu} /X > t f(x))}{ P(X > x)}
\rightarrow 0 
$$
since $f(x)  \to \infty$.
Assumption \ref{firstcony} is verified similarly.
For Assumption \ref{secondcon}, choose $L =1$ and as $x \to \infty$,
$$
\frac{P( X > f(x), Y > f(x))}{ P( X  > x)} = \frac{P( X >
f(x), e^{2\mu} /X > f(x))}{ P(X > x, Y > x)}
\rightarrow  0.
$$
We conclude  by Theorem \ref{mainthm},
\[ P( X + Y > x) \sim 2P(X > x) .\] 
}
\end{ex}

\begin{ex} \label{heavylight} 
{\rm Example \ref{rhoextreme}
is a special case of a more general phenomenon. Suppose, $F \in MDA(\Lambda)$ with auxiliary function $f(x)$ having the property 
\begin{equation} \label{exeqn4}
\liminf_{ x\to \infty} f(x) =\delta > 0. 
\end{equation}
Assume that the support of $F$ is a subset of $[0, \infty)$ and $x_0 = \sup\{ x : F(x) < 1\} = \infty$ and also that
 $x_1 = \inf\{ x : F(x) >  0 \} = 0$. Distributions satisfying these conditions include the exponential,
gamma, lognormal.
Define
$ X = F^{\leftarrow} (U) ,$ and $Y = F^{\leftarrow} (1 - U)$, where $U \sim \rm{Uniform}(0,1)$.
This pair $(X, Y)$ satisfies the Assumptions in Section
\ref{ass_sec}. 

Checking Assumption \ref{ratioxy} is easy since $X$ and $Y$ are
identically distributed. 
To verify Assumption \ref{firstconx}, fix $t  > 0$ and
define $\epsilon_t = F( \frac{t \delta}{2} )$. Since, $x_1 = 0$, we
have $\epsilon_t > 0$. Then, for large $x$ making $f(x)>\delta/2$, we have
\begin{align*}
&\frac{P( X > x, Y > t f(x))}{ P(X > x)} = \frac{P( U  > F(x), 1 - U > F(t f(x)))}{ P(X > x)}\\
&\qquad\le \frac{P( U  > F(x), 1 - U > \epsilon_t)}{ P(X > x)} 
= \frac{P( U  > F(x), U < 1- \epsilon_t)}{ P(X > x)}
\rightarrow 0 
\end{align*}
since $F(x)  \to 1$, and $ x_0 = \infty$.
Assumption \ref{firstcony} is similarly verified. To verify Assumption
\ref{secondcon}, 
choose $L$ such that $F( \frac{ L \delta}{2} ) > \frac{1}{2}$ and 
for $x$ sufficiently large, 
\begin{align*}
&\frac{P( X > L f(x), Y > L f(x))}{ P( X  > x)} \le \frac{P( X >
\frac{ L \delta}{2}, Y > \frac{ L \delta}{2})}{ P(X > x)} \\
&\qquad = \frac{P( U >F( \frac{ L \delta}{2}), 1- U > F(\frac{ L \delta}{2})
  )}{ P(X > x)}= 0. 
\end{align*}
Hence, $(X, Y)$ satisfy the Assumptions of Section \ref{ass_sec} and
by Theorem \ref{mainthm}, 
\[ P( X + Y > x) \sim 2P(X > x). \] 

In this example, if $\lim_{x \to \infty} f(x) = \infty$, we do not
need the condition $x_1 = 0$.
}
\end{ex}

\begin{rem}
{\rm Note, in the previous two examples a comonotonic dependence structure is used. }
\end{rem}

\begin{ex} \label{lognormal} 
{\rm Suppose,
$ X = \exp(X_1) ,$ $ Y = \exp(X_2), $
where $(X_1, X_2)$ is bivariate normal with correlation $\rho \in [-1,1)$.
For simplicity, assume  each $X_i$ has mean $\mu$ and variance
$\sigma^2 > 0 $. 
This example is extensively considered in
\cite{asmussen:nandayapa}. We have already considered the case $ \rho =
-1$ in Example \ref{rhoextreme}, so here we consider $\rho
\in ( -1, 1)$.

Assumptions \ref{Gumbel} and \ref{ratioxy} of Section \ref{ass_sec}
are easily verified. Following the same reason as in
Example \ref{rhoextreme}, we take the auxiliary function to be  
$$f(x) = \frac{ \sigma^2 x}{ \log x - \mu}.$$ 
Observe,
\begin{align}\label{exlogneqn1}
\frac{\log f(x) - \mu}{ \sigma} =& \frac{ \log  \left(\frac{ \sigma^2 x}{ \log x - \mu}\right) - \mu}{\sigma}
= \frac{\log x - \mu}{ \sigma}  - \frac{1}{\sigma}\log  \left(\frac{ \log x - \mu}{ \sigma^2} \right) \nonumber \\
=& \left( \frac{\log x - \mu}{ \sigma} \right)( 1 + o(1)).
\end{align}
For Assumption \ref{firstconx} , we have for  $t > 0$, as $x \to \infty$,
\begin{align*}
&\frac{P( X > x, Y > t f(x))}{ P(X > x)} = \frac{P( X_1 > \log
x, X_2  > \log t f(x))}{ P(X_1 > \log x)}\\
&\qquad \le \frac{P( X_1 + X_2 > \log x + \log (t f(x)))}{ P(X_1 > \log x)}\\
&\qquad = \frac{\bar \Phi \left( \frac{1}{\sqrt{2 \sigma^2 (1+ \rho)}}( \log
    x + \log (t 
 f(x)) - 2\mu) \right)}{ \bar \Phi( \frac{\log x -\mu}{\sigma})} \\
&\qquad = \frac{\bar \Phi \left( \frac{1}{\sqrt{2 (1+ \rho)}}( \frac{\log x -\mu}{\sigma} + \frac{\log f(x) -\mu}{\sigma} + \frac{\log t}{\sigma} )\right)}{ \bar \Phi( \frac{\log x -\mu}{\sigma})} \\ 
&\qquad = \frac{\bar \Phi \left( \frac{2}{\sqrt{2 (1+ \rho)}}( \frac{\log x
      -\mu}{\sigma} )( 1 + o(1))\right)}{ \bar \Phi( \frac{\log x
    -\mu}{\sigma})} 
\rightarrow 0,
\end{align*}
where we used \eqref{exlogneqn1} and  the fact that $\Phi \in
MDA(\Lambda)$ and therefore $\bar \Phi$ is
$-\infty$-varying    \citep[page 53]{resnick:1987}. Note, 
$\rho < 1$ entails $\frac{2}{\sqrt{2(1 + \rho)}} > 1$.

For Assumption \ref{secondcon}, choose $L = 1$. As $x \to \infty$, we
have using  \eqref{exlogneqn1},
\begin{align*}
&\frac{P( X > f(x), Y > f(x))}{ P(X > x)} = \frac{P( X_1  >
\log f(x), X_2 > \log f (x))}{ P(X_1 > \log x )}\\
&\qquad \le \frac{P( X_1 + X_2 > 2 \log f(x))}{ P(X_1  > \log x)}
 = \frac{\bar \Phi \left( \frac{2( \log f(x) - \mu)}{\sqrt{2\sigma^2(1+ \rho)}} 
 \right)}{ \bar \Phi (\frac{(\log x - \mu)}{\sigma})} \\
&\qquad = \frac{\bar \Phi \left( \frac{2}{\sqrt{2 (1+ \rho)}}(
    \frac{\log x -\mu}{\sigma} )( 1 + o(1))\right)}{ \bar \Phi(
  \frac{\log x -\mu}{\sigma})} 
 \rightarrow 0.
 \end{align*}
 }
\end{ex}

\begin{ex}
 {\rm Let $X_1$ and $X_2$ be independent and identically
  distributed with the common distribution $H \in MDA(\Lambda)$,
  having auxiliary function $f_1(\cdot)$ satisfying
  \eqref{exeqn4} and infinite right end point. 
Also, suppose, $F \in MDA( \Lambda)$ with auxiliary function $f_2(
\cdot)$, concentrates on 
$[0,\infty)$ and satisfies the conditions in Example \ref{heavylight}.
Define $X$ and $Y$ as
\[ X = F^{\leftarrow}(U) \wedge X_1, \hskip 0.5 cm Y = F^{\leftarrow}(1-U) \wedge X_2 \]
where $U$ is a uniformly distributed random variable on (0,1)
which is independent of $(X_1, X_2)$.

{}From Proposition 1.4 of \cite[page 43] {resnick:1987}, the
distribution of $X$  belongs to the maximal domain of attraction of
Gumbel with auxiliary function  
\[ f(x) = \frac{f_1(x)f_2(x)}{ f_1(x) + f_2(x) } \]
Hence,
\[ \limsup_{x \to \infty} \frac{1}{f(x)} \le  \limsup_{x \to \infty} \frac{1}{f_1(x)} + \limsup_{x \to \infty} \frac{1}{f_2(x)}  < \infty, \]
and thus,
\[\liminf_{x \to \infty} f(x) > 0. \]
Also, note,
\begin{align*}
P( X > x ) =& P( U > F(x) , X_1 > x) = P( U > F(x))P( X_1 > x)
 =\bar F (x) \bar H(x) \\
\intertext{and}
 P( Y > x ) =& P( 1 - U > F(x) , X_2 > x) = P( 1 - U > F(x))P( X_2 >
 x)
 = \bar F (x) \bar H(x).
\end{align*}
Arguing as in Example \ref{heavylight}, we can show that the pair $(X,
Y)$ satisfy  the assumptions in Section \ref{ass_sec}.
}  
\end{ex}

\begin{ex}
{\rm Here is an example of a distribution for $(X, Y)$ where our
  assumptions are not satisfied, but the asymptotic behavior is the
  same as in Theorem \ref{mainthm}. Suppose, $X$ and $Y$ are iid with
  common distribution $F$, where 
\[ \bar F (x) = \exp( -x^{\alpha} ) \hskip 0.5 cm \alpha \in (0,1), \hskip 0.5 cm x > 0. \]
This distribution is extensively studied in \cite{rootzen:1986} and
satisfies $F \in MDA( \Lambda) \cap \mathcal{S}$. Since it is subexponential,
\[ P( X + Y > x) \sim 2P( X > x). \]
However, this distribution does not satisfy Assumption \ref{secondcon}
of Section \ref{ass_sec}. 

Since $F$ is a Von-Mises function, we may take the auxiliary function to be
\[ f(x) = \frac{ \bar F (x)}{ F^{'}(x)}  = \frac{x^{ 1 - \alpha}}{ \alpha}.  \]
 Assumption \ref{secondcon} is not satisfied for any $L >
0$, since for any $L  > 0$, as $x \to \infty$,
\begin{align*}
&\frac{P(X > L f(x), Y > L f(x) )}{ P( X > x)} = \frac{{\left[\bar F(L
      f(x)) \right]}^2}{ \bar F (x)} 
=  \frac{ \exp( -2{[L f(x)]}^{\alpha}) }{  \exp( -x^{\alpha} )}\\
&\qquad =  \frac{ \exp( -2{(\frac{L}{\alpha})}^{\alpha} x^{ \alpha(1 -
    \alpha)}) }{  \exp( -x^{\alpha} )}
= \exp\left( x^{\alpha}( 1- 2{(\frac{L}{\alpha})}^{\alpha} x^{ -
    \alpha^2})\right)
\rightarrow \infty .
\end{align*} 
 This also shows the criteria \eqref{subeqn1} for $F \in
\mathcal{S}$ is sufficient but not necessary.
}   
\end{ex}

\section{Linear combinations of random variables with non-negative coefficients}
This section studies linear combinations of risks $X, Y$ with non-negative coefficients. We consider two cases:
(i)  the distributions of $X$ and $Y$ are tail-equivalent,
and (ii) the distributions of $X$ and $Y$ lack tail-equivalence.
We explicitly give the asymptotic tail behavior of the linear combinations
of risks in the tail-equivalent case and also in
one special case where  tail-equivalence is absent. We note
 that one cannot expect similar behavior in the two cases. 

\subsection{Tail-equivalent cases}
\subsubsection{Linear combination of two random variables with non-negative coefficients}
 \begin{thm} \label{lincombthm} 
Assume, $(U, V)$ is a pair of random
variables which satisfy Assumptions \ref{Gumbel}, \ref{firstconx}, \ref{firstcony} and \ref{secondcon} of Section \ref{ass_sec}.
Moreover, assume that Assumption \ref{ratioxy} holds in the form
\begin{equation} \label{linthmeqn0}
\lim_ {x \rightarrow \infty} \frac{ P( V > x )}{ P( U > x)} = c
\in (0, \infty).
\end{equation}
Define, $\hat S_2 = a_1U + a_2V$ and $ a_i \ge 0, i = 1, 2$ and set $m_2 = a_1 \vee a_2$. Then, as $ x \to \infty$,
 \[ P( \hat S_2 > x ) \sim P( U > \frac{x}{m_2} )\left[1_{\{a_1 = m_2\}} + c 1_{\{a_2 = m_2\}} \right]  \]
\end{thm}

We assume $U$ and $V$ are tail equivalent, i.e. the
constant $c$ cannot be 0 and hence
both the marginal distributions belong to $MDA(\Lambda)$, the maximal
domain of attraction of the Gumbel.
If  $\lim_{x \to 
\infty} {P(V > x)}/{P(U > x)} = 0$, the asymptotic behavior of 
$P( a_1U + a_2V > x)$ as $x\to \infty$ can be different as
 illustrated in the following example.

\begin{ex}
{\rm Assume, $(U, V)$ are iid random
variables with common distribution $F$, which satisfy Assumptions
\ref{Gumbel}, \ref{firstconx}, \ref{firstcony} and
\ref{secondcon} of Section \ref{ass_sec}. Define the  two 
random vectors by $(U_1, V_1) = ( U, \frac{1}{5}V)$ and $(U_2,
V_2) = ( U, \frac{1}{2}V)$.  Both pairs
$(U_1, V_1)$ and  $(U_2, V_2)$ satisfy Assumptions \ref{Gumbel},
 \ref{firstconx}, \ref{firstcony} and \ref{secondcon}
of Section \ref{ass_sec}. For both pairs, c = 0,
i.e., 
$$
\lim_ {x \rightarrow \infty} \frac{ P( V_1 > x )}{ P( U_1 > x)} = 0
\quad \text{ and } \quad
\lim_ {x \rightarrow \infty} \frac{ P( V_2 > x )}{ P( U_2 > x)} = 0.
$$
Since, $(U, V)$ satisfies the Assumptions of Theorem \ref{lincombthm}, we have as $x \to \infty$,
\begin{align*}
 P ( 3U_1 + 10V_1 > x ) &=  P( 3U + 2V >  x ) \sim  P( 3U > x ) =  P(
 3U_1  > x ), \\
\intertext{and}
 P (3U_2 + 10V_2 > x ) &=  P( 3U + 5V >  x ) \sim  P( 5V > x ) =  P(
 10V_2  > x ) .
 \end{align*}
This example illustrates we cannot expect Theorem \ref{lincombthm} to hold for the case $c = 0$.
}
\end{ex}

We now turn to the proof of Theorem \ref{lincombthm}.

\begin{proof}
The case $a_1
= a_2$  is resolved by Theorem
 \ref{mainthm} since
$$
 P ( a_1(U + V) > x ) =  P( U + V > \frac{x}{a_1} ) 
\sim (1 + c)P ( U >  \frac{x}{a_1} ).
$$
So the interesting cases are $a_1 >  a_2 $ and $ a_1 <  a_2$ and for
the following, assume $a_1 > a_2$, the  other case being similar.

 There is nothing to prove if $a_2 = 0$, so assume $a_1 > a_2 >
  0$ which makes $ m_2 = a_1$. It suffices to check the Assumptions in Section \ref{ass_sec} for $X = U$
  and $Y = a_2 V/a_1$. For this definition of $X$, $Y$, we have
  \begin{equation}\label{linthmeqn1}
\lim_{x \to \infty}\frac{P(Y > x)}{P(X > x)} = \lim_{x \to
  \infty}\frac{P(a_2V/a_1 > x)}{P(U > x)}
                        = \lim_{x \to \infty}\frac{P( V >
                          a_1x/a_2)}{P(U > x)} =
                         0. 
  \end{equation}
  The last equality is true from \eqref{linthmeqn0} and the fact
that the tail of any distribution in $MDA(\Lambda)$
 is $-\infty$-varying \citep[page 53]{resnick:1987}. 
{}From Theorem \ref{mainthm} and \eqref{linthmeqn1}, we  get, as $x \to
\infty$, 
\begin{align*}
P( a_1U +& a_2V > x) = P( a_1(U + a_2V/a_1) >x)\\
=& P( U + a_2V/a_1 > x/a_1 )
= P( X + Y > x/a_1 )\\
\sim & P( U > x/a_1) 
= P( U > \frac{x}{m_2} )\left[1_{\{a_1 = m_2\}} + c 1_{\{a_2 =
m_2\}} \right].
\end{align*}

To complete the proof, the
Assumptions in Section \ref{ass_sec} must be verified for $X = U$ and
$Y = a_2V/a_1$. 
Assumption \ref{Gumbel} is assumed in the Theorem and 
Assumption \ref{ratioxy} was verified in \eqref{linthmeqn1}. 
For Assumption \ref{firstconx}, note that $U \in MDA(\Lambda)$
and suppose $f(\cdot)$ is the auxiliary function of the
 distribution of $U$. By hypothesis, for   $t > 0$,
 \begin{equation} \label{linthmeqn2}
\lim_{ x \to \infty} P( |V| > t f(x) | U > x ) = 0,
\end{equation}
and therefore, using \eqref{linthmeqn2},
$$
\lim_{ x \to \infty}P( a_2|V|/a_1 > t f(x) | U > x ) = \lim_{ x \to
  \infty}P( |V| > a_1 t f(x)/a_2 | U > x)
= 0.
$$

  Remark \ref{rem:asy_ind}(\ref{c=0}) implies we do not need to verify 
Assumption \ref{firstcony}, so we check
Assumption \ref{secondcon}. For this we have, as $x\to \infty$,
\begin{align*}
\frac{P( a_2V/a_1 > L f(x),U > L f(x) )}{P(U > x)} &= \frac{P( V >
a_1L f(x)/a_2, U > L f(x) )}{P(U >x)}\\
&\le \frac{P( V > L f(x), U > L f(x) )}{P(U >x)}
\rightarrow 0.       
\end{align*}
This proves the case  $ a_1 > a_2$.
\end{proof}

\subsubsection{Linear Combination of more than two random variables with non-negative coefficients}
 
 \begin{cor} \label{lincor1}
Assume, $X_1, X_2, \ldots X_d$ are non-negative random
variables which pairwise satisfy Assumptions \ref{firstconx},
\ref{firstcony}, \ref{secondcon} of Section \ref{ass_sec}. 
Further suppose the distribution of $X_1$ satisfies Assumption
\ref{Gumbel} of Section \ref{ass_sec} and that
\begin{equation}\label{lincor1eqn0}
\lim_ {x \rightarrow \infty} \frac{ P( X_i > x )}{ P( X_1 > x)} =
c_i \in (0, \infty), \hskip 0.5 cm i = 1,2,\ldots,d.
\end{equation}
Set $c_1 =1$ and define for $d>1$, $\hat S_d = a_1X_1 + a_2X_2 + \ldots
a_dX_d,$ for  $ a_i \ge 0, \hskip 0.2 cm i= 1, 2, \ldots, d$. Also, define, 
$$
m_d = \bigvee_{i=1}^d a_i  \quad  \text{and}\quad 
 N_d = \sum_{\{ 1 \le i \le d : a_i = m_d \}} c_i.
$$
Then
 \[ P( \hat S_d > x ) \sim N_d P( X_1 > \frac{x}{m_d} ),\qquad x\to\infty. \]
\end{cor}

This result is consistent with the case where $X_1, X_2, \ldots, X_d$ are
iid with common distribution in
$MDA( \Lambda) \cap \mathcal{S} $; see \cite{davis:resnick:1988}.

The random variables  $X_1, X_2, \ldots X_d$ 
are tail-equivalent and satisfy Assumption \ref{firstconx} of Section
\ref{ass_sec} pairwise. Therefore Remark \ref{rem:asy_ind}(\ref{asy_ind})
implies pairwise asymptotic independence and hence, by \cite[page
291]{resnick:1987}, $X_1\dots,X_d$ are asymptotically 
independent.

In the special case that the random variables are identically
distributed, 
 $N_d = |\{ 1 \le i \le d : a_i = m_d \}|$, where $| \cdot |$
is the size of a set.

\begin{rem}
\rm{ It is possible to prove Corollary \ref{lincor1} using Corollary \ref{sumcor}. However, in the proof it is usually difficult to verify Assumption \ref{firstcony} of Section \ref{ass_sec}. Note,  a similar problem is avoided carefully in the proof of Theorem \ref{lincombthm} through the help of Remark  \ref{rem:asy_ind}(\ref{c=0}). Though a remark similar to Remark \ref{rem:asy_ind}(\ref{c=0}) could  also be made for Corollary \ref{sumcor}, it is notationally inconvenient. So, to avoid this notational difficulty, Theorem \ref{lincombthm} is used for the proof.}
\end{rem}
\begin{proof}
 Proceeding by induction, note the base case for $d = 2$ is proved in Theorem \ref{lincombthm}. As an induction hypothesis, suppose the result is true for $d = k$, so, as $x \to \infty$,
$$
 P( \hat S_k > x) \sim N_k P( X_1 > \frac{x}{m_k} ) 
\sim \frac{N_k}{c_{k + 1}} P( X_{k + 1} > \frac{x}{m_k} ).
$$
To prove the result for $d = k + 1$, notice,
 \begin{eqnarray} \label{lincor1eqn1}
 m_{k + 1} = m_k \vee a_{k+1}, 
 \end{eqnarray}
 and
 \begin{align}\label{lincor1eqn2}
 N_{k +1} =  \left[c_{k + 1}1_{\{a_{k+1} = m_{k+1} \}} + N_k 1_{\{m_k = m_{k+1}\}}
 \right] , \nonumber\\
 \intertext{ so that}
  \frac{N_{k +1}}{c_{k +1}} =  \left[1_{\{a_{k+1} = m_{k+1} \}} + \frac{N_k}{c_{k + 1}} 1_{\{m_k = m_{k+1}\}} \right].
 \end{align}
 By the induction hypothesis,
\begin{align}\label{lincor1eqn3}
\lim_{ x\to \infty} \frac{ P( m_k^{-1}\hat S_k >  x )}{ P (X_{k + 1}
> x) } = \lim_{ x\to \infty} \frac{ P( m_k^{-1}\hat S_k >  x )}{ P (X_1
> x) } \frac{ P (X_1 > x) }{ P (X_{k + 1} > x) } 
= \frac{N_k}{c_{k + 1}}. 
\end{align}
If we prove the assumptions in  Theorem \ref{lincombthm}  are valid
with $U = X_{k +1}$ and $V = m_k^{-1} \hat S_k$, then, Theorem
\ref{lincombthm}, \eqref{lincor1eqn1},  \eqref{lincor1eqn2} and
\eqref{lincor1eqn3} imply, as $x \to \infty$, 
\begin{align*}
P(\hat S_{k+1}& > x) = P( a_{k+1}X_{k+1} + m_k \hat S_k > x)\\
&\sim \frac{N_{k+1}}{c_{k + 1}}P( X_{k + 1} > \frac{x}{m_{k+1}})
\sim N_{k+1}P( X_1 > \frac{x}{m_{k + 1}}), 
\end{align*}
and by induction, our result holds for all $d \ge 2$.

Assumption \ref{Gumbel} is assumed.
For \eqref{linthmeqn0}, consider that on the one hand,
\begin{align}  \label{lincor1eqn4}
 N_k &= \sum_{\{ 1 \le i \le k : a_i = m_k \}} c_i 
\, \ge \bigwedge_{i = 1}^k c_i  > 0
\end{align}
and on the other,
\begin{align} \label{lincor1eqn5}
 N_k = \sum_{\{ 1 \le i \le k : a_i = m_k \}} c_i 
 \le k \bigvee_{i = 1}^k c_i < \infty,
\end{align}
and therefore the limit in \eqref{lincor1eqn3}
satisfies $N_k/c_{k+1} \in (0, \infty).$

 Next, suppose, two random variables $U$ and $V$ are tail equivalent
 and both belong to $MDA(\Lambda)$.
 If $f(\cdot)$, $\tilde f(\cdot)$ are the auxiliary
 functions of $U$ and $V$ respectively, then 
$f(x) \sim \tilde f(x)$, as $x \to \infty$; see
\cite{resnick:1971, resnick:1971b} . Since, in the present case, all the
random variables are 
tail-equivalent, 
Remark \ref{rem:asy_ind}(\ref{mean_excess}) 
implies we  can work with the
auxiliary function of any one of them, say $X_{k +1}$. So, 
 $X_1, X_2, \ldots, X_d$  satisfy 
Assumptions \ref{firstconx}, \ref{firstcony} and \ref{secondcon} of
Section \ref{ass_sec} pairwise with the auxiliary function $f( \cdot )$ of 
$X_{k + 1}$. That is, for $1 \le i \ne j \le d$, and any $t >
0$,
\begin{equation} \label{lincor1eqn6}
\lim_{x \to \infty} P( X_j > t f(x) | X_i > x )= 0
\end{equation}
and for $1 \le i < j \le d$, for some $L_{ij} > 0$
\begin{equation} \label{lincor1eqn7}
\frac{P( X_i > L_{ij} f(x), X_j > L_{ij} f(x) )}{P(X_i
>x)} = 0.
\end{equation}

To verify Assumption \ref{firstconx}, observe for 
$t  >  0$, that as  $x \to \infty $,
\begin{align*}
     P( &|m_k^{-1}\hat S_k| > t f(x) | X_{ k + 1} > x )\\
&\le P(a_1X_1 + a_2X_2 + \dots + a_kX_k > m_kt  f(x) | X_{ k + 1} > x )\\
   &\le \sum_{ i = 1}^{k}  P( X_i  > a_i^{-1}m_k t f(x)/k | X_{ k + 1} > x )
\le \sum_{ i = 1}^{k}  P( X_i  > t f(x)/k | X_{ k + 1} > x )
 \rightarrow 0,
\end{align*}
by \eqref{lincor1eqn6}.
For Assumption \ref{firstcony}, note,
\begin{eqnarray}\label{lincor1eqn8}
\lim_{x \to \infty} \frac{P( m_k^{-1}\hat S_k > x)}{P(X_1 > x)} &=&
N_k,
\end{eqnarray}
and for $1 \le i \le k$,
\begin{eqnarray}\label{lincor1eqn9}
\lim_{x \to \infty}\frac{P( (m_k^{-1}\hat S_k > x) \cap
(m_k^{-1}a_iX_i > x))}{P(X_1 > x)} &=& \lim_{x \to \infty}\frac{ P(m_k^{-1}a_iX_i > x)}{P(X_1 > x)} \nonumber \\
&=& c_i 1_{\{a_i = m_k\}}.
\end{eqnarray}
 The first equality uses the assumption that $X_i$'s are
non-negative. The second equality is true from \eqref{lincor1eqn0} and the fact
that the tail of any distribution in the maximal domain of
attraction of Gumbel is $-\infty$-varying. Now, for
$1 \le i < j \le k$, using \eqref{asymp_ind},
\begin{align*}
 \lim_{x \to \infty}&\frac{P( (m_k^{-1}\hat S_k > x) \cap (
m_k^{-1}a_iX_i > x) \cap (m_k^{-1}a_jX_j > x))}{P(X_1 > x)}
 \\ &\le \lim_{x \to \infty} \frac{P(( m_k^{-1}a_iX_i
> x) \cap (m_k^{-1}a_jX_j > x))}{P(X_1 > x)}  \\
 &\le \lim_{x \to \infty} \frac{ P((X_i > x) \cap (X_j > x)) }{P(X_1 > x)}  = 0.
\end{align*}
Therefore, using \eqref{lincor1eqn9}, 
\begin{align}\label{lincor1eqn10}
 \lim_{ x \to \infty}& \frac{P( (m_k^{-1}\hat S_k > x) \cap
(\cup_{i=1}^k(m_k^{-1}a_iX_i
> x)))}{P(X_1 > x)}  \nonumber \\   &= \lim_{ x \to \infty} \frac{\sum_{i=1}^k P((m_k^{-1}\hat S_k >
x)\cap (m_k^{-1}a_iX_i > x))}{P(X_1 > x)}= N_k .
\end{align}
{}From \eqref{lincor1eqn8} and \eqref{lincor1eqn10} it follows that
$$
\lim_{x \to \infty}\frac{P( (m_k^{-1}\hat S_k > x) \cap
{(\cup_{i=1}^k(m_k^{-1}a_iX_i
> x))}^c)}{P(X_1 > x)} = 0,
$$
and this, along with \eqref{lincor1eqn0} and \eqref{lincor1eqn3} give
\begin{equation}\label{lincor1eqn11}
\lim_{x \to \infty}\frac{P( (m_k^{-1}\hat S_k > x) \cap
{(\cup_{i=1}^k(m_k^{-1}a_iX_i
> x))}^c)}{P(m_k^{-1}\hat S_k > x)} = 0.
\end{equation}
Now, we check Assumption \ref{firstcony}. For
$t  >  0$, as $x \to \infty$,
\begin{align*}
    P( &|X_{ k +1}| > t f(x) | m_k^{-1} \hat S_k > x ) = \frac{P( X_{
        k + 1 } > t f(x), m_k^{-1}\hat S_k > x)}{ P( m_k^{-1}\hat S_k
      > x )}\\ 
  & \sim  \frac{P( X_{ k + 1 } > t f(x), m_k^{-1}\hat S_k > x, \cup_{i=1}^k \{ m_k^{-1}a_iX_i > x\} )}
   { P( m_k^{-1}\hat S_k > x )} \\
   &\le \frac{P( X_{ k + 1 }|> t f(x), \cup_{i=1}^k \{ m_k^{-1}a_iX_i > x\} )}{ P( \hat S_k > m_k x )}\\
&\le \frac{\sum_{i=1}^k P( X_{ k + 1 } > t f(x), m_k^{-1}a_iX_i > x )}{ P( \hat S_k > m_k x
   )},\\
   \intertext{where we have used \eqref{lincor1eqn11}. Using our
     induction hypothesis, we get that the quantity
     above is aymptotically equivalent to} 
   &\sim \frac{\sum_{i=1}^k P( X_{ k + 1 } > t f(x), m_k^{-1}a_iX_i > x )}{ N_kP( X_1 > x )} \\
    &\le \frac{\sum_{i=1}^k P( X_{ k + 1 } > t f(x), X_i > x )}{ N_kP(
      X_1 > x )} \\
& = \frac{\sum_{i=1}^k P( X_{ k + 1 } > t f(x), X_i > x )}{P( X_i > x )}\frac{P( X_i > x )}{ N_kP( X_1 > x
    )}
      \rightarrow 0,
\end{align*}
by  \eqref{lincor1eqn6}.

For assumption \ref{secondcon}, let, $L = kL_{max}$, where $L_{max} =
\max_{1 \le i \le k}L_{i, k + 1}$ (recall, equation \eqref{lincor1eqn7}) . Then  using
\eqref{lincor1eqn3}, \eqref{lincor1eqn11} and \eqref{lincor1eqn7}, we have
 \begin{align*}
 & \frac{P( X_{ k +1} > kL_{max} f(x) , m_k^{-1} \hat S_k > kL_{max} f(x) )}{ P(X_{k+1} > x )} \\
&\qquad    \sim  \frac{P( X_{ k + 1 } > kL_{max} f(x), m_k^{-1}\hat S_k > x, \cup_{i=1}^k \{ m_k^{-1}a_iX_i > x\} )}
   {P( X_{k+1} > x )} \\
&\qquad   \le \frac{P( X_{ k + 1 } > k L_{max} f(x), \cup_{i=1}^k \{ m_k^{-1}a_iX_i > L_{max}f(x) \})}{ P( X_{k+1} > x )}\\
&\qquad   \le \frac{\sum_{ i = 1}^k P( X_{ k + 1 } > L_{i, k + 1}f(x), m_k^{-1}a_iX_i > L_{i, k + 1}f(x) )}{ P( X_{k+1} > x )} \\
&\qquad   \le \ \frac{\sum_{ i = 1}^k P( X_{ k + 1 } > L_{i, k + 1}f(x),
  X_i > L_{i, k + 1}f(x) )}{ P( X_{k+1} > x )}  \rightarrow  0.
   \end{align*}
\end{proof}

\subsection{One special case where the distributions are possibly NOT tail-equivalent}

\begin{thm} \label{splcase} 
Assume, $Y_1, Y_2, \ldots Y_d$ are identically distributed
non-negative random variables. Also, assume $ a_i, 
\beta_i \ge 0,$    $ i= 1, 2, \ldots, 
d$. For $d \ge 1$, define, $\hat S_d = a_1Y_1^{\beta_1}
+ a_2Y_2^{\beta_2} + \ldots + a_dY_d^{\beta_d}$ and set 
\[ \beta = \bigvee_{i=1}^d \beta_i, \hskip 1 cm  q_d = \bigvee_{\{ 1
\le i \le d : \beta_i = \beta \}} a_i, \]\[ J_d =|\{ 1 \le i \le d :
\beta_i = \beta, a_i = q_d \}| \] where $|\cdot|$ denotes the size of
the set. Suppose, $q_dY_1^{\beta}, q_dY_2^{\beta}, \ldots
q_dY_d^{\beta}$ pairwise satisfy Assumptions \ref{firstconx},
\ref{firstcony} and \ref{secondcon} of Section \ref{ass_sec} and that the
distribution of $q_dY_1^{\beta}$ satisfies Assumption \ref{Gumbel}
of Section \ref{ass_sec} where the auxiliary function $f(x)$ satisfies the additional condition that $f(x) \to \infty$, as $x \to \infty$. Then,
 \[ P( \hat S_d > x ) \sim J_d P( Y_1^{\beta} > \frac{x}{q_d} ) \]
\end{thm}

\begin{rem} \label{not-tail-equivalent} {\rm If $\beta_1 > \beta_2$, then
$Y_1^{\beta_1}$ and $Y_2^{\beta_2}$ are not tail-equivalent. Note, in
this case, the asymptotic approximation of $P( a_1Y_1^{\beta_1} +
a_2Y_2^{\beta_2} > x)$ does not depend on $a_2$.

Theorem \ref{splcase} shows  different tail behavior from the
tail-equivalent cases but follows the paradigm that only the heaviest
tails matter. The Theorem  shows that Theorem 1 of
\cite{asmussen:nandayapa}  is a special case of a more general
phenomenon. Let  $(X_1, X_2, \ldots, X_d) \sim 
N( \bf{0}, \Sigma ) $, where  
\[ \Sigma = ( \rho_{ij}) , \hskip 0.2 cm \rho_{ii} = 1,
\hskip 0.1 cm  \forall \hskip 0.1 cm i , \hskip 0.1 cm \rho_{ij} <
\hskip 0.1 cm   1 \le i <  j \le d . \] 
Let, $(Y_1, Y_2, \ldots, Y_d) \sim (\exp(X_1), \exp (X_2), \ldots, \exp (X_d))$. Clearly, 
\[ a_iY_i^{\beta_i} \sim \rm{Lognormal} ( \log a_i, \beta_i^2) \] 
From Example \ref{lognormal}, $(q_dY_1^{\beta}, q_dY_2^{\beta}, \ldots, q_dY_d^{\beta})$ satisfy the Assumptions of Theorem \ref{splcase}, where $q_d, \beta$ have the same meaning as in Theorem \ref{splcase}. Also, $( Z_1, Z_2, \ldots, Z_d) = ( a_1Y_1^{\beta_1}, a_2Y_2^{\beta_2}, \ldots, a_nY_d^{\beta_d})$ satisfies the Assumptions of Theorem 1 of \cite{asmussen:nandayapa}. The results of that theorem and Theorem \ref{splcase} match.} 
\end{rem}

\begin{proof} Without loss of generality, assume $\beta_1 = \beta$ and $a_1
= q_d$. Also, assume $a_i > 0$ for $i
= 1,2, \ldots, d$. Denote,
\[ X_i = a_iY_i^{\beta_i}\hskip 1cm i = 1, 2, \ldots, d. \]
To start, suppose, for some $i \in \{ 2, \ldots, d\}, \beta_i <
\beta$. Then, for large  $x$, $[ a_iY_i^{\beta_i} > x ] \subseteq [ \frac{q_d}{2}Y_i^{\beta} > x ]$, and hence for large  $x$,
\begin{align*}
P(a_iY_i^{\beta_i} > x ) \le P(\frac{q_d}{2}Y_i^{\beta} > x) = P(q_dY_1^{\beta} > 2x). 
\end{align*}

Then,
\begin{equation}\label{linthm2eqn1}
c_i = \lim_{x \to \infty}\frac{P(X_i > x)}{P(X_1 > x)} 
= \lim_{x \to \infty}\frac{P(a_iY_i^{\beta_i} > x)}{P(q_dY_1^{\beta} > x )} 
\le \lim_{x \to \infty}\frac{P(q_dY_1^{\beta} > 2x)}{P(q_dY_1^{\beta} > x)} 
                       = 0.
  \end{equation}

Next, suppose, for some $i \in \{ 2, \ldots, d\}, \beta_i = \beta, a_i < q_d$. Then,
\begin{equation}\label{linthm2eqn2}
c_i = \lim_{x \to \infty}\frac{P(X_i > x)}{P(X_1 > x)}
= \lim_{x \to \infty}\frac{P(a_iY_i^{\beta} > x)}{P(q_dY_1^{\beta} > x)} = \lim_{x \to
  \infty}\frac{P(q_dY_1^{\beta} > \frac{q_dx}{a_i})}{P(q_dY_1^{\beta} > x)} = 0. 
  \end{equation}
  In both the equations \eqref{linthm2eqn1} and \eqref{linthm2eqn2}, the last equalities are true from the fact
that the tail of any distribution in the maximal domain of
attraction of the Gumbel is $-\infty$-varying.

Finally, suppose, for some $i \in \{ 2, \ldots, d\}, \beta_i = \beta, a_i = q_d$.
\begin{equation}\label{linthm2eqn3}
c_i = \lim_{x \to \infty}\frac{P(X_i > x)}{P(X_1 > x)} = \lim_{x \to
  \infty}\frac{P(Y_i^{\beta} > \frac{x}{q_d})}{P(Y_1^{\beta} >
  \frac{x}{q_d})} = \lim_{x \to \infty}\frac{P(Y_1^{\beta} >
  \frac{x}{q_d})}{P(Y_1^{\beta} > \frac{x}{q_d})} =1. 
 \end{equation}

  It suffices to check the assumptions in Corollary \ref{sumcor} with
  this set of $X_1, X_2, \ldots, 
  X_d$, since then Corollary \ref{sumcor} and  \eqref{linthm2eqn1},
  \eqref{linthm2eqn2}, \eqref{linthm2eqn3} would imply, as $x \to
  \infty$, 
$$
P( \hat S_d > x) \sim (1 + \sum_{i=2}^d c_i)P( X_1 > x) \sim J_d P(
X_1 > x ) = J_d P( Y_1^{\beta} > \frac{x}{q_d} ) .
$$

Assumption \ref{Gumbel} is assumed in the Theorem statement
and \eqref{cor1eqn1} is
already shown in \eqref{linthm2eqn1}, \eqref{linthm2eqn2} and
\eqref{linthm2eqn3}. For assumptions \ref{firstconx} and
\ref{firstcony}, proceed as follows. 
 By hypothesis, we know that $X_1$ belongs to the maximal domain of attraction of the Gumbel distribution.
 Let $f(\cdot)$ be the auxiliary function corresponding to the
 distribution of $X_1$. By hypothesis, we know, for $t > 0$, for $1 \le i \ne  j \le d$, 
 \begin{equation}
\lim_{ x \to \infty} P( q_dY_j^{\beta} > t f(x) | q_dY_i^{\beta} > x
)= 0.
\end{equation}
Using Remark \ref{asyindpair}(\ref{relax}), it is enough to show
\begin{equation} \label{linthm2eqn4}
\lim_{ x \to \infty} \frac{P( X_j > t f(x) , X_i > x)}{P(X_1
> x)} = 0 ,
\end{equation}
and to see this, note that since $f(x) \to \infty$, for large $x$ and for all $t > 0$, $[ X_j > t f(x) , X_i > x ]\subseteq [ q_dY_j^{\beta} > t f(x), q_dY_i^{\beta} > x ]$. Hence,
\begin{align*}
\lim_{ x \to \infty}  \frac{P( X_j > t f(x) , X_i > x)}{P(X_1
> x)} \le  \lim_{ x \to \infty} \frac{P( q_dY_j^{\beta} > t f(x) ,
q_dY_i^{\beta} > x)}{P(q_dY_i^{\beta}
> x)} = 0.
\end{align*}

For Assumption \ref{secondcon}, using Remark
\ref{asyindpair}(\ref{relax}), 
we show, for some $ L >0$,
\begin{align}  \label{linthm2eqn5}
\lim_{ x \to \infty} \frac{P( X_j > L f(x) , X_i > L f(x))}{P(X_1
> x)} &= 0.
\end{align}
By hypothesis, we know, for all \hskip 0.1 cm $1 \le i < j \le d$, there exists some
$L_{ij} > 0$,
\begin{equation}  \label{linthm2eqn6}
\lim_{ x \to \infty} \frac{P( q_dY_j^{\beta} > L_{ij} f(x) ,
q_dY_i^{\beta}
> L_{ij} f(x) )}{P( q_dY_i^{\beta} > x)} = 0.
\end{equation}
Also, note that since $f(x) \to \infty$, for large $x$,  $[ X_j > L_{ij} f(x), X_i > L_{ij}
f(x) ] \subseteq [ q_dY_j^{\beta} > L_{ij} f(x) , q_dY_i^{\beta} > L_{ij} f(x) ] $. Hence,
\begin{align*}
\lim_{ x \to \infty} &\frac{P( X_j > L_{ij} f(x), X_i > L_{ij}
f(x))}{P(X_1 > x)} \le  \lim_{ x \to \infty} \frac{P( q_dY_j^{\beta} > L_{ij} f(x) ,
q_dY_i^{\beta} > L_{ij} f(x))}{P(q_dY_1^{\beta}
> x)}   =  0 
\end{align*}
by  \eqref{linthm2eqn6}.
\end{proof}

\section{An Optimization Problem}

\subsection{The problem}

Suppose, we have a portfolio consisting of $d$ financial
instruments. The risk per unit of the $i$-th instrument is $X_i$. The
goal is to earn revenue  \$L. Assume, each unit of the $i$-th instrument
earns $\$l_i$ over the chosen time horizon. Subject to 
earnings being at least \$L, how many units of each instrument,
$a_1, a_2, \ldots, a_d$, should be used to build the
portfolio, so that the probability that the total portfolio risk 
$a_1X_1 + a_2X_2 + \ldots + a_dX_d$ exceeds some fixed
large threshold $x$, is minimal?

Thus, consider the following optimization problem:
\begin{align*}
 \min _{\{a_1,\dots,a_d\}} & P\left[ \sum_{i=1}^d a_iX_i >  x \right]  \\
 \text{s.t.}\quad & a_1l_1 + a_2l_2 + \dots
+ a_dl_d \ge L, \\
{}& a_i \ge 0, \; i = 1,2, \ldots, d. 
\end{align*}

For a more general case, consider the following optimization
problem:
\begin{align*}
 \min _{\{a_1,\dots,a_d\}}& P\left[ \sum_{i=1}^d a_iX_i >  x \right] \\
 \text{s.t.}\quad & h(a_1,a_2, \ldots, a_d) \ge L, \\
{}& a_i \ge 0, \; i = 1,2, \ldots, d. 
\end{align*}

\subsection{The method}

Suppose, $X_1, X_2, \ldots, X_d$ satisfy the assumptions of
Corollary \ref{lincor1}. Even with these assumptions, exact solution
of the optimization problem is difficult. An obvious way to obtain an
approximate solution to the optimization problem is
to assume that the threshold $ x $ is big and use the
 asymptotic approximation of $P(a_1X_1 + a_2X_2 + \ldots
+ a_dX_d > x)$ from Corollary \ref{lincor1}, hoping that 
the solution of the resulting optimization problem is close to the actual optimal
value.
So, using the notation of Corollary \ref{lincor1}, we solve the
following optimization problem: 
\begin{align*}
 \min _{\{a_1,\dots,a_d\}} & N_d P( X_1 > \frac{x}{m_d} )  \\
\text{s.t.}\quad & h(a_1,a_2, \ldots, a_d) \ge L, \\
{}& a_i \ge 0, \; i = 1,2, \ldots, d. 
\end{align*}

Suppose $ \hat a_1, \hat a_2, \ldots , \hat a_d$ and 
$\tilde a_1, \tilde a_2,
\ldots, \tilde a_d$ 
are two feasible solutions  for the
given set of constraints. Set
\[ \hat m_d =  \bigvee_{i=1}^d \hat a_i , \hskip 1cm  \hat N_d = \sum_{ \{ 1 \le i \le d : \hat a_i = \hat m_d \}} c_i \]
\[ \tilde m_d =  \bigvee_{i=1}^d \tilde a_i , \hskip 1cm  \tilde N_d = \sum_{ \{ 1 \le i \le d : \tilde a_i = \tilde m_d \}} c_i \]
If $\hat m_d > \tilde m_d$, then since, $P[X_1\leq x] \in
MDA(\Lambda)$, as $x \to \infty$, 
\[ \frac{P( X_1 > x/ \hat m_d )}{ P( X_1 > x/ \tilde m_d )}
\rightarrow \infty .\] 
Now, since  both $\hat N_d, \tilde N_d \in [ \wedge_{i=1}^d c_i, d
\vee_{i=1}^d c_i ] $, we have as $x \to \infty$, 
\[ \frac{ \hat N_dP( X_1 > x/ \hat m_d )}{\tilde N_d P( X_1 > x/ \tilde m_d )}  \rightarrow \infty. \]
So, we hope that
$\tilde a_1, \tilde a_2, \ldots, \tilde a_d$ is a better 
feasible solution for  the optimization problem.

Thus, values of $a_1,a_2, \ldots, a_d$ which solve the
above optimization problem can be computed by solving the following
two optimization problems in sequence.
\begin{enumerate}
\item[(i)] First solve
\begin{align*}
\min _{\{a_1\dots,a_d\} }  m_d &= \max\{ a_1,a_2, \ldots, a_d \}  \\
 \text{s.t.}\qquad & h(a_1,a_2, \ldots, a_d) \ge L, \\
{}& a_i \ge 0, \;i = 1,2, \ldots, d.
\end{align*}

\item[(ii)] Suppose, the best choice of $a_1, a_2, \ldots, a_d$ gives m as the value
of the objective function
for the  optimization problem in (i).
Then we solve
\begin{align*}
 \min _{\{a_1,\dots,a_d\}}  N_d &= \sum_{\{ 1 \le i \le d : a_i = m \}} c_i  \\
 \text{s.t.}\qquad &h(a_1,a_2, \ldots, a_d) \ge L, \\
{}& \max\{a_1,a_2, \ldots, a_d\} = m, \\
{}& a_i \ge 0\; i = 1,2, \ldots, d. 
\end{align*}
\end{enumerate}

\subsection{ A special case}

The motivating case is that  $h$ is a linear function with positive
coefficients
of the form
\[ h(a_1,a_2, \ldots, a_d) = a_1l_1 + a_2l_2 + \ldots + a_dl_d .\]
The approximate solution using the asymptotic form of $P[\sum_{i=1}^d
a_iX_i>x]$ is
\[ a_1 = a_2 = \ldots = a_d = L/(l_1 +l_2 + \ldots + l_d) .\]
This leads to $m = L/(l_1 +l_2 + \ldots + l_d)$ and $N_d = \sum_{i=1}^d c_i$.

\section{Simulation studies}\label{sec:sim}
We carried out some simulation studies to check for fixed large
thresholds the accuracy of the
asymptotic approximation in Theorem \ref{mainthm} and also
to check how good is the approximate solution  for
the optimization problem.  As  expected,  in some
cases the approximation works  well whereas in others it
performs  poorly which suggests caution about using the
asymptotic results for numerical purposes.
Simulation also suggests that the
approximate solution of the optimization problem works well
in cases where simulation studies
suggest that the approximation is good for fixed large thresholds. One
particular model studied,  Example \ref{lognormal} with $\mu = 0,
\sigma =1$, is noted here to illustrate the point. We varied
$\rho$ and observed the asymptotic behavior of the sum of the risks.

\subsection{Where is the approximation good?}\label{subsec:good}
To test the approximation for $P( X + Y > x)$,
we have to find good simulation estimates of the probabilities  $P( X + Y > x)$.
This, however is not easy, especially in the case
when the marginal distributions of the risks $X$ and $Y$ are subexponential and is
still a topic of current research in the simulation community. The
approach usually taken in these 
cases is Conditional Monte Carlo \citep[page 173]{asmussen:glynn:2007}. So, this method is used to compute 
$P( X + Y > x)$ and the simulation estimates  are compared with the theoretical
approximations.

 The simulation of $P( X + Y > x)$ uses the algorithm
 suggested in \cite{asmussen:nandayapa} for $\rho \in (-1,1)$ who also 
 note the   properties of this algorithm.
 If $\rho = -1$, we have a way to compute the 
probability exactly.  In this case, $X = 1/Y$ almost surely, so
in the following manner we compute the required probability:
\begin{align*}
P\Bigl(& X + \frac{1}{X}  > x  \Bigr) = P\Bigl( X > \frac{x + \sqrt{x^2-2}}{2}\Bigr)
+ P\Bigl( X < \frac{x - \sqrt{x^2-2}}{2} \Bigr)\\
=&  P\Bigl( \log  X > \log {\Bigl(\frac{x + \sqrt{x^2-2}}{2}\Bigr)}\Bigr) + P\Bigl(  \log X <
\log{\Bigl(\frac{x -
\sqrt{x^2-2}}{2}\Bigr)}\Bigr)\\
= & \bar \Phi \Bigl(\log {\Bigl(\frac{x + \sqrt{x^2-2}}{2}\Bigr)}\Bigr) + \Phi
\Bigl(\log {\Bigl(\frac{x - \sqrt{x^2-2}}{2}\Bigr)}\Bigr)
\end{align*}

\subsubsection{Patterns in the results}\label{subsec:patterns}
 For judging the quality of the asymptotic approximation, we
focus on the simulation estimate $P( X + Y > x )$ and not the
threshold $x$, since a change of  distribution may imply a change in
how rare is a particular threshold crossing.
So, when comparing the quality
of the asymptotic approximation across different models, it makes more sense to focus
on the value of $P( X +Y > x)$, rather than the particular threshold $x$.
When $\rho = -1$, exact calculations
suggest that the approximation is extremely good even when
the actual probability $P( X + Y > x)$ is of the order of
$10^{-2}$. As expected, the asymptotic approximation 
improves as a function of increasing  threshold.
When $\rho \in (-1,1)$, we rely on the
simulation estimate as a surrogate for the exact tail probability and
compare it with the theoretical approximations.

The results indicate that the closer $\rho$  is  to $-1$, the better
the approximation.
For $\rho =-1$, the approximation is good for events with probability of
the order of $10^{-2}$ and to achieve comparable  precision in the
relative error when $\rho =0$, the event has to be much rarer and 
have a probability of the order of $10^{-10}$. For $\rho
= 0.9$, the results for different thresholds did not show any
convergence pattern. This emphasizes that in practice the numerical 
approximations should be
used with caution. Clearly for $\rho = 1$ the asymptotic approximation
is not correct and  $\rho = 0.9$ is expected to behave somewhat
like the case when $\rho =1$.

The tables give representative results.
We first give the results for  $\rho = -1$ in Table \ref{table:1}, since in this case
no simulation is required. The column `Ratio' in Table \ref{table:1} is defined as
$$
\text{Ratio} =  \frac{\hbox{Actual probability}}{\hbox{
Asymptotic approximation}}.
$$

\begin{table}[t]
\caption{$\rho=-1$}
\begin{center}
\label{table:1}
\begin{tabular}{|c|c|c|c|}
\hline Threshold & Actual probability & Asymptotic approximation & Ratio\\
\hline 10 & 0.0219 & 0.0213 & 1.0272\\
\hline 16 & 0.0056 & 0.0056 & 1.0121\\
\hline 24 & 0.0015 & 0.0015 & 1.0060\\
\hline 30 & 6.7365 $\times 10^{-4}$ & 6.7091 $\times 10^{-4}$ &1.0041\\
\hline 100 & 4.1233 $\times 10^{-6}$ & 4.1213 $\times 10^{-6}$ &
1.0005\\
\hline 1000 & 4.9238 $\times 10^{-12}$ & 4.9238 $\times 10^{-12}$ &
1.0000\\
\hline
\end{tabular}
\end{center}
\end{table}

For subsequent tables, the columns `Ratio' and `Half-width' are defined as 
\begin{align*}
\text{Ratio} = &\frac{\hbox{Simulation estimated probability}}{\hbox{ Asymptotic approximation}}\\
\text{Half-width} = & \text{Half-width of the 95\% confidence interval
  of the ratio.}
\end{align*}
 In each case, $10^7$ observations were used to compute the probability estimates.

\begin{table}[h]
\label{table:-.9}
\caption{$\rho=-0.9$}
\begin{center}
\begin{tabular}{|c|c|c|c|c|c|}
\hline Threshold & Simulation estimated probability & Asymptotic approximation &
Ratio &  Half-width \\
\hline 3 & 0.3687 & 0.2719 & 1.3556 & 0.0006\\
\hline 5 & 0.1207 & 0.1075 & 1.1227 &  0.0012\\
\hline 10 & 0.0221 & 0.0213 & 1.0375 & 0.0026\\
\hline 20 & 0.0028 & 0.0027 & 1.0082 & 0.0064\\
\hline 30 & 6.8873 $\times 10^{-4}$ & 6.7091 $\times 10^{-4}$ &
1.0265 &  0.0119\\
\hline 40 & 2.2134 $\times 10^{-4}$ & 2.2524 $\times 10^{-4}$ &
0.9827 &  0.0183\\
\hline 50 & 9.3675 $\times 10^{-5}$ & 9.1526 $\times 10^{-5}$ &
1.0235 &  0.0285\\
\hline
\end{tabular}
\end{center}
\end{table}

\begin{table}[h]
\label{table:0}
\caption{$\rho = 0$}
\begin{center}
\begin{tabular}{|c|c|c|c|c|c|}
\hline Threshold & Simulation estimated probability & Asymptotic approximation &
Ratio &  Half-width \\
\hline 10 & 0.0338 & 0.0213 & 1.5844 & 0.0033\\
\hline 50 & 1.0798 $\times 10^{-4}$ & 9.1526 $\times 10^{-5}$ & 1.1798 & 0.0002\\
\hline 100 & 4.5032 $\times 10^{-6}$ & 4.1213 $\times 10^{-6}$ &
1.0927  & 0.0001\\
\hline 300 & 1.2117 $\times 10^{-8}$ & 1.1718 $\times 10^{-8}$ &
1.0341 & 0.0000\\
\hline 600 & 1.6147 $\times 10^{-10}$ & 1.5853 $\times 10^{-10}$ &
1.0185  & 0.0122\\
\hline 1000 & 4.9821 $\times 10^{-12}$ & 4.9238 $\times 10^{-12}$ &
1.0118 &  0.0000\\
\hline 2000 & 1.9620 $\times 10^{-14}$ & 2.9310 $\times 10^{-14}$ &
1.0106 & 0.0000\\
\hline
\end{tabular}
\end{center}
\end{table}

\begin{table}
\label{table:.9}
\caption{$\rho = 0.9$}
\begin{center}
\begin{tabular}{|c|c|c|c|c|c|}
\hline Threshold & Simulation estimated probability & Asymptotic approximation &
Ratio &  Half-width \\
\hline 10 & 0.0521 & 0.0213 & 2.4439 & 
0.0088\\
\hline 30 & 0.0030 & 6.7091 $\times 10^{-4}$ & 4.4081 & 0.0275\\
\hline 50 & 5.2652 $\times 10^{-4}$ & 9.1526 $\times 10^{-5}$ &
5.7527 &  0.0759\\
\hline 75 & 1.1217 $\times 10^{-4}$ & 1.5781 $\times 10^{-5}$ &
7.1077 &  0.1843\\
\hline 100 & 3.4333 $\times 10^{-5}$ & 4.1213 $\times 10^{-6}$ &
8.3307 & 0.3642\\
\hline
\end{tabular}
\end{center}
\end{table}

\subsection{ How good is the portfolio suggestion?}
Here, we consider the quality of our approximate solutions for the
optimization problem. We choose the same risk model given in Example \ref{lognormal},
because we have information about which values of $\rho$ lead to
good  asymptotic
approximation.  We resort to a naive
method for analyzing the optimization.
 For different $(a_1,a_2)$, we obtain
estimates of $P( a_1X + a_2Y > x)$ through simulation. 
To get the estimates proceed as follows: For 
 $a_1, a_2 > 0$
\[ \left(\begin{array}{c}a_1X \\ a_2Y \end{array} \right) =
\left(\begin{array}{c}\exp\{ \log(a_1) + X_1 \}\\ \exp\{ \log(a_2) + X_2\}
\end{array} \right)\]
Now,
\[ \left(\begin{array}{c} Z_1 \\  Z_2
\end{array} \right) = \left(\begin{array}{c} \log(a_1) + X_1 \\  \log(a_2) + X_2
\end{array} \right) \sim \mathcal{N} \left( \left(\begin {array}{c} \log(a_1)\\\log(a_2) \end {array}\right), \left(\begin {array}{cc} 1 & \rho \\ \rho & 1 \end {array}\right)\right) \hskip 0.5 cm \rho \in [-1, 1) \]
So, again we are in the framework of  \cite{asmussen:nandayapa}, and
we use the algorithm given in their paper to estimate the rare event probabilities. 
When either $a_1$ or $a_2$ is equal to 0, we can
compute the exact probability and hence do not need an estimate. We choose
$(a_1, a_2)$ in the following way. Let $C$ be the set of all possible
$(a_1, a_2)$ which satisfy the constraint. First, $a_1$ is chosen from
the corresponding projection of $C$ with a small grid, and then for
each $a_1$, $a_2$ is determined from the constraint. Let us call this
set $C^*$. For $(a_1, a_2) \in C^*$, $P( a_1X + a_2Y > x)$ is
estimated through simulation and then it is observed which $(a_1,
a_2)$ gives the minimum estimate of $P( a_1X + a_2Y > x)$. Let, $( \tilde a_1, \tilde a_2)$ be this pair; i.e. $P( \tilde a_1 X + \tilde a_2 Y > x)  = \min_{ ( a_1, a_2) \in C^*} P( a_1 X + a_2 Y > x) $. Also, let
$(a_1^*, a_2^*)$ be the approximate solution of the optimization
problem as noted in the previous section. Relative error of the
approximate solution is computed by comparing $P( a_1^* X + a_2^* Y >
x)$ with $\min_{ (a_1, a_2) \in C^*} P( a_1X + a_2Y > x)$. 

\subsubsection{Identifying  patterns}
We do not have  error estimates for our
simulation results. One could consider  bootstrapping to obtain such error estimates,
but we have not done so.  Despite the weaknesses
of the naive procedure, the results  are interesting.

We note one case with the  linear constraint $2a_1 + 3a_2 = 1$. 
The suggested optimum  portfolio based on asymptotic approximation
is $(a_1^*, a_2^*) = (0.2, 0.2)$. The 3 cases where
$\rho = -0.9, 0, 0.9$, are chosen, the reason being
that we know from the results in earlier section that the asymptotic
approximation is good in the case $\rho = -0.9$, reasonable when
$\rho = 0$ and rather bad when $\rho = 0.9$.
The approximate solution $(a_1^*, a_2^*)$
relies on replacing the original objective function by its
asymptotic approximation, and so it is reasonable to
expect different accuracies for these three values of $\rho$
and this turned out to be the case.
In the cases of $\rho = -0.9$ and $\rho= 0$, we
see that $\tilde a_1$ comes close to $0.2$ as the
threshold $x$
increases. But, in the case of $\rho = 0.9$, no pattern in the convergence of $\tilde a_1$ is observed which is expected because for  $\rho =1$,
both the risks are actually the same random variable which implies
indifference to the choice of $(a_1, a_2) \in C$.

Another remark is that in each  case 
of $\rho = -0.9, 0, 0.9$, the relative
errors do not show any convergence pattern. Perhaps to expect
otherwise is unrealistic
as we are using the minimum of some simualtion
estimates to compute the relative error. Still, 
we illustrate through an example the accuracy  by
comparing with an extreme case where we build the portfolio
consisting of only  one  asset. For $\rho = 0$, and  threshold
$x=10$, the extreme cases will yield probabilities $0.2441$ and
$0.1360$. These risk probabilities are quite high compared that of our
suggested optimal portfolio $(a_1^*, a_2^*)$ based on asymptotic
approximation, which has risk probability $P(  a_1^* X + a_2^* Y > x)
= 1.0793 \times 10^{-4}$;  also, the minimum of the simulation
estimates $P( \tilde a_1 X + \tilde a_2 Y > x)$ is of the same
order. So, the suggested portfolio $(a_1^*, a_2^*)$ is quite effective
in reducing the risk and possibly close to the best one.

The following additional conclusion can be made.
In the case of $\rho =
-0.9$, even when $P( \tilde a_1 X + \tilde a_2 Y > x)$ is as big as
$0.11$, it is quite close  to $P( a_1^* X +  a_2^* Y > x)$, indicating that the suggested optimal choice $(a_1^*, a_2^*)$ significantly reduces the risk  probability. 
For $\rho = 0$, a comparable statement can be made 
when the minimum of the probability
estimates is of the order of $10^{-2}$. However,
for $\rho = 0.9$, the relative errors are never  close to 0. 
Interestingly, even for $\rho = 0.9$, $P( \tilde a_1 X + \tilde a_2 Y > x)$ and $P( a_1^* X +  a_2^* Y > x)$ are almost always of the same order. However, it should be noted at this point that even
in this case of $\rho = 0.9$, the extreme cases where the portfolio is built on entirely one of the assets, $P( a_1 X +  a_2 Y > x)$ is of a much bigger order than $P( \tilde a_1 X + \tilde a_2 Y > x)$. So, in this case, possibly $P( a_1 X +  a_2 Y > x)$ differs considerably  from choices where $a_1, a_2 > 0$ and the case where either $a_1 =0$ or $a_2 = 0$, but does not differ too much among the choices where $(a_1, a_2) \in C, a_1, a_2 > 0$.  This fact  justifies the intuition as mentioned before that the case $\rho = 0.9$ is similar to case $\rho =1$. Some of the results are
noted in tables below.

Results are summarized in the tables for 
 $\rho = -0.9, 0, 0.9$ and
 constraint $2a_1 + 3a_2 = 1$. For each fixed $\rho$, we give
\begin{itemize}
\item the threshold $x$, 
\item $\tilde a_1$, where $( \tilde a_1, \tilde a_2)
  \in C^*$ and
$$P( \tilde a_1X + \tilde a_2Y > x) = \min_{ (a_1, a_2)
    \in C^*} P( a_1X + a_2Y > x),$$
\item  $E1= \min_{ (a_1, a_2) \in C^*} P(  a_1X + a_2Y > x)$, 
\item  $E2 = P( a_1^* X + a_2^* Y > x)$, 
\item  the `Relative error' = $\frac{E2 -E1}{E1}$.
\end{itemize}
For each value of $\rho$,  $a_1$ is chosen
with gap 0.01 from the projection of $C^*$; i.e. we considered ($a_1 =0,
0.01, 0.02, \ldots$, 0.5). For
each such $a_1$, we used 10000 observations to obtain the estimates of
the probability $P( a_1X + a_2Y > x)$.

\begin{table}[h]
\label{tableopt:-0.9}
\caption{$\rho = -0.9$}
\begin{center}
\begin{tabular}{|c|c|c|c|c|}
\hline Threshold & $\tilde a_1$  & E1 &
E2 & Relative error \\
\hline 1 & 0.13 & 0.1097 & 0.1204 & 0.0975\\
\hline 3 & 0.18 & 0.0067 & 0.0069 & 0.0322 \\
\hline 5 & 0.19 & 0.0013& 0.0013 & 0.0294 \\
\hline 10 & 0.19 & 1.0299 $\times 10^{-4}$ & 1.0592 $\times 10^{-4}$ & 0.0284\\
\hline 20 & 0.21 & 2.0806 $\times 10^{-6}$ & 2.0806 $\times 10^{-6}$ & 1.2213 $\times 10^{-15}$\\
\hline
\end{tabular}
\end{center}
\end{table}

\begin{table}
\label{tableopt:0}
\caption{$\rho = 0$}
\begin{center}
\begin{tabular}{|c|c|c|c|c|}
\hline Threshold & $\tilde a_1$  & E1 &
E2 & Relative error \\
\hline 1 & 0.03 & 0.1349 & 0.1723 & 0.2765\\
\hline 3 & 0.16 & 0.0093 & 0.0101 & 0.0759 \\
\hline 5 & 0.18 & 0.0016 & 0.0017 & 0.0503 \\
\hline 10 & 0.19 & 1.0424 $\times 10^{-4}$ & 1.0793 $\times 10^{-4}$ & 0.0354\\
\hline 20 & 0.20 & 4.3888 $\times 10^{-6}$ & 4.3888 $\times 10^{-6}$ & 0\\
\hline
\end{tabular}
\end{center}
\end{table}

\begin{table}
\label{tableopt:0.9}
\caption{$\rho = 0.9$}
\begin{center}
\begin{tabular}{|c|c|c|c|c|}
\hline Threshold & $\tilde a_1$  & E1 &
E2 & Relative error \\
\hline 1 & 0.01 & 0.1360 & 0.1798 & 0.3223\\
\hline 3 & 0.01 & 0.0140 & 0.0208 & 0.4831 \\
\hline 5 & 0.02 & 0.0033 & 0.0050 & 0.5146 \\
\hline 10 & 0.02 & 2.8357 $\times 10^{-4}$ & 4.9475 $\times 10^{-4}$ & 0.7447\\
\hline 20 & 0.04 & 1.3241 $\times 10^{-6}$ & 2.4023 $\times 10^{-6}$ & 0.8142\\
\hline
\end{tabular}
\end{center}
\end{table}

\section{Concluding Remarks}
An important case for the study of asymptotic behavior of the sum of
 risks is the case where the risks are asymptotically independent,
identically distributed and belong to the maximal domain of attraction
of the Gumbel distribution. Many commonly occuring risk 
distributions fall in this category.
We have provided  sufficient conditions for
\[ \lim_{ x \to \infty} \frac{ P( X + Y > x)}{P(X > x)} = 2,\]
and extended the conditions
to cover the case where the marginal distributions are not the same
and to the case where some risk distributions  have lighter tail but
the distribution does not belong to the 
maximal domain of attraction of the Gumbel. We are not able to provide
necessary and sufficient conditions for this kind of asymptotic
behavior which is an unresolved problem.
It will be interesting to see if it is possible to find a
distribution of risks $(X,Y)$ for which  the risks are asymptotically
independent, identically distributed, belong to $MDA(\Lambda)$,
 and the asymptotic behavior of the sum is
different than two cases mentioned in the introduction, viz.  
\[ \lim_{ x \to \infty} \frac{ P( X + Y > x)}{P(X > x)} \in \{2, \infty\} \] 

Even for cases where the asymptotic behavior is understood,
nothing is  known about the rate of convergence in these cases; i.e. a
 quantitative estimate how good the approximation $2P( X > x)$ is for
 the quantity $P( X + Y > x)$ for a large threshold $x$. 
Simulation studies indicate in certain circumstances the 
 approximation is accurate, but in other cases its accuracy is dismal.

We observed in the previous section that when
tail probability approximation is good, the approximate solution of
the optimization problem is also accurate
 whereas in the other cases
this solution has poor accuracy.
 So, results on the rate of convergence
would contribute to understanding the appropriateness of the approximate
 solutions in different scenarios.

 An anonymous and conscientious referee provided many insightful and helpful comments.

\bibliographystyle{plainnat}
\bibliography{bibfile}
  \end{document}